\documentclass[a4paper]{amsart}
\usepackage{color}
\usepackage[latin1]{inputenc}
\usepackage[T1]{fontenc}
\usepackage{amsfonts}
\usepackage{amssymb}
\usepackage{amsmath}
\usepackage{amsthm}
\usepackage{stmaryrd}
\usepackage{eufrak}
\usepackage{graphicx,type1cm,eso-pic,color}
\usepackage{pstricks,pst-plot,pstricks-add}
\usepackage{float}
\newtheorem{theorem}{Theorem}[section]
\newtheorem{lemma}[theorem]{Lemma}
\newtheorem{proposition}[theorem]{Proposition}
\newtheorem{corollary}[theorem]{Corollary}

\newtheorem{assumption}[theorem]{Assumption}
\newtheorem{remark}[theorem]{Remark}
\begin{document}
\setlength\arraycolsep{2pt}
\title{Kernel Selection in Nonparametric Regression}
\author{H\'el\`ene HALCONRUY*}
\address{*LTCI, T\'el\'ecom Paris, Palaiseau, France}
\email{helene.halconruy@telecom-paris.fr}
\author{Nicolas MARIE$^\dag$}
\address{$^\dag$Laboratoire Modal'X, Universit\'e Paris Nanterre, Nanterre, France}
\email{nmarie@parisnanterre.fr}
\address{*,$^\dag$ESME Sudria, Paris, France}
\email{helene.halconruy@esme.fr}
\email{nicolas.marie@esme.fr}
\keywords{Nonparametric estimators ; Projection estimators ; Model selection ; Regression model.}
\date{}
\maketitle
\noindent
%


%
\begin{abstract}
In the regression model $Y = b(X) +\sigma(X)\varepsilon$, where $X$ has a density $f$, this paper deals with an oracle inequality for an estimator of $bf$, involving a kernel in the sense of Lerasle et al. (2016), selected via the PCO method. In addition to the bandwidth selection for kernel-based estimators already studied in Lacour, Massart and Rivoirard (2017) and Comte and Marie (2020), the dimension selection for anisotropic projection estimators of $f$ and $bf$ is covered.
\end{abstract}
\tableofcontents
\noindent
\textbf{MSC2010:} 62G05 ; 62G08.
%


%
\section{Introduction}\label{section_introduction}
Consider $n\in\mathbb N^*$ independent $\mathbb R^d\times\mathbb R$-valued ($d\in\mathbb N^*$) random variables $(X_1,Y_1),\dots,(X_n,Y_n)$, having the same probability distribution assumed to be absolutely continuous with respect to Lebesgue's measure, and
\begin{displaymath}
\widehat s_{K,\ell}(n;x) :=
\frac{1}{n}\sum_{i = 1}^{n}K(X_i,x)\ell(Y_i)
\textrm{ $;$ }x\in\mathbb R^d,
\end{displaymath}
where $\ell :\mathbb R\rightarrow\mathbb R$ is a Borel function and $K$ is a symmetric continuous map from $\mathbb R^d\times\mathbb R^d$ into $\mathbb R$. This is an estimator of the function $s :\mathbb R^d\rightarrow\mathbb R$ defined by
\begin{displaymath}
s(x) :=
\mathbb E(\ell(Y_1)|X_1 = x)f(x)
\textrm{ $;$ }
\forall x\in\mathbb R^d,
\end{displaymath}
where $f$ is a density of $X_1$. For $\ell = 1$, $\widehat s_{K,\ell}(n;.)$ coincides with the estimator of $f$ studied in Lerasle et al. \cite{LMRB16}, covering Parzen-Rosenblatt's and projection estimators already deeply studied in the literature (see Parzen \cite{PARZEN62}, Rosenblatt \cite{ROSENBLATT56}, Tsybakov \cite{TSYBAKOV09}, etc.), but for $\ell\not= 1$, it covers estimators involved in nonparametric regression. Assume that for every $i\in\{1,\dots,n\}$,
\begin{equation}\label{nonparametric_regression}
Y_i = b(X_i) +\sigma(X_i)\varepsilon_i
\end{equation}
where $\varepsilon_i$ is a centered random variable of variance $1$, independent of $X_i$, and $b,\sigma :\mathbb R^d\rightarrow\mathbb R$ are Borel functions.
\begin{itemize}
 \item If $\ell =\textrm{Id}_{\mathbb R}$, $k$ is a symmetric kernel and
 \begin{equation}\label{PR_kernel}
 K(x',x) =\prod_{q = 1}^{d}\frac{1}{h_q}k\left(\frac{x_q' - x_q}{h_q}\right)
 \textrm{ with }
 h_1,\dots,h_d > 0
 \end{equation}
 for every $x,x'\in\mathbb R^d$, then $\widehat s_{K,\ell}(n;.)$ is the numerator of the well-known Nadaraya-Watson estimator of the regression function $b$ (see Nadaraya \cite{NADARAYA64} and Watson \cite{WATSON64}). Precisely, $\widehat s_{K,\ell}(n;.)$ is an estimator of $s = bf$ because $\varepsilon_1$ is independent to $X_1$ and $\mathbb E(\varepsilon_1) = 0$. If $\ell\not=\textrm{Id}_{\mathbb R}$, then $\widehat s_{K,\ell}(n;.)$ is the numerator of the estimator studied in Einmahl and Mason \cite{EM00,EM05}.
 \item If $\ell =\textrm{Id}_{\mathbb R}$, $\mathcal B_{m_q} =\{\varphi_{1}^{m_q},\dots,\varphi_{m_q}^{m_q}\}$ ($m_q\in\mathbb N^*$ and $q\in\{1,\dots,d\}$) is an orthonormal family of $\mathbb L^2(\mathbb R)$ and
 \begin{equation}\label{projection_kernel}
 K(x',x) =
 \prod_{q = 1}^{d}
 \sum_{j = 1}^{m_q}\varphi_{j}^{m_q}(x_q)\varphi_{j}^{m_q}(x_q')
 \end{equation}
 for every $x,x'\in\mathbb R^d$, then $\widehat s_{K,\ell}(n;.)$ is the projection estimator on $\mathcal S =\textrm{span}(\mathcal B_{m_1}\otimes\dots\otimes\mathcal B_{m_d})$ of $s = bf$.
\end{itemize}
Now, assume that $b = 0$ in Model (\ref{nonparametric_regression}): for every $i\in\{1,\dots,n\}$,
\begin{equation}\label{heteroscedastic_model}
Y_i =\sigma(X_i)\varepsilon_i.
\end{equation}
If $\ell(x) = x^2$ for every $x\in\mathbb R$, then $\widehat s_{K,\ell}(n;.)$ is an estimator of $s =\sigma^2f$.
\\
\\
These ten last years, several data-driven procedures have been proposed in order to select the bandwidth of Parzen-Rosenblatt's estimator ($\ell = 1$ and $K$ defined by (\ref{PR_kernel})). First, Goldenshluger-Lepski's method, introduced in \cite{GL11}, which reaches the adequate bias-variance compromise, but is not completely satisfactory on the numerical side (see Comte and Rebafka \cite{CR16}). More recently, in \cite{LMR17}, Lacour, Massart and Rivoirard proposed the PCO (Penalized Comparison to Overfitting) method and proved an oracle inequality for the associated adaptive Parzen-Rosenblatt's estimator by using a concentration inequality for the U-statistics due to Houdr\'e and Reynaud-Bouret \cite{HRB03}. Together with Varet, they established the numerical efficiency of the PCO method in Varet et al. \cite{LMRV}. Still in the density estimation framework, the PCO method has been extended to bandwidths selection for the recursive Wolverton-Wagner estimator in Comte and Marie \cite{CM20}.
\\
Comte and Marie \cite{CM} deal with an oracle inequality and numerical experiments for an adaptive Nadaraya-Watson's estimator with a numerator and a denominator having distinct bandwidths, both selected via the PCO method. Since the output variable in a regression model has no reason to be bounded, there were significant additional difficulties, bypassed in \cite{CM}, to establish an oracle inequality for the numerator's adaptive estimator. Via similar arguments, the present article deals with an oracle inequality for $\widehat s_{\widehat K,\ell}(n;.)$, where $\widehat K$ is selected via the PCO method in the spirit of Lerasle et al. \cite{LMRB16}. As in Comte and Marie \cite{CM}, one can deduce an oracle inequality for the adaptive quotient estimator $\widehat s_{\widehat K,\ell}(n;.)/\widehat s_{\widehat L,1}(n;.)$ of $\mathbb E(\ell(Y_1)|X_1 =\cdot)$, where $\widehat K$ and $\widehat L$ are both selected via the PCO method.\\
In addition to the bandwidth selection for kernel-based estimators already studied in \cite{LMR17,CM}, the present paper covers the dimension selection for projection estimators of $f$, $bf$ when $Y_1,\dots,Y_n$ are defined by Model (\ref{nonparametric_regression}) with $\ell =\textrm{Id}_{\mathbb R}$, and $\sigma^2f$ when $Y_1,\dots,Y_n$ are defined by Model (\ref{heteroscedastic_model}) with $\ell(x) = x^2$ for every $x\in\mathbb R$. For projection estimators, when $d = 1$, the usual model selection method (see Comte \cite{COMTE14}, Chapter 2, Section 5) seems hard to beat. However, when $d > 1$ and $K$ is defined by (\ref{projection_kernel}), $m_1,\dots,m_d$ are selected via a Goldenshluger-Lepski type method (see Chagny \cite{CHAGNY13}), which has the same numerical weakness than the Goldenshluger-Lepski method for bandwidth selection when $K$ is defined by (\ref{PR_kernel}). So, for the dimension selection for anisotropic projection estimators, the PCO method is interesting.
\\
\\
In Section \ref{section_risk_bound}, some examples of kernels sets are provided and a risk bound on $\widehat s_{K,\ell}(n;.)$ is established. Section \ref{section_kernel_selection} deals with an oracle inequality for $\widehat s_{\widehat K,\ell}(n;.)$, where $\widehat K$ is selected via the PCO method.
%


%
\section{Risk bound}\label{section_risk_bound}
Throughout the paper, $s\in\mathbb L^2(\mathbb R^d)$. Let $\mathcal K_n$ be a set of symmetric continuous maps from $\mathbb R^d\times\mathbb R^d$ into $\mathbb R$, of cardinality less or equal than $n$, fulfilling the following assumption.
%


%
\begin{assumption}\label{assumption_K}
There exists a deterministic constant $\mathfrak m_{\mathcal K,\ell} > 0$, not depending on $n$, such that
\begin{enumerate}
 \item For every $K\in\mathcal K_n$,
 \begin{displaymath}
 \sup_{x'\in\mathbb R^d}
 \|K(x',.)\|_{2}^{2}
 \leqslant
 \mathfrak m_{\mathcal K,\ell}n.
 \end{displaymath}
 \item For every $K\in\mathcal K_n$,
 \begin{displaymath}
 \|s_{K,\ell}\|_{2}^{2}\leqslant\mathfrak m_{\mathcal K,\ell}
 \end{displaymath}
 with
 \begin{displaymath}
 s_{K,\ell} :=\mathbb E(\widehat s_{K,\ell}(n;.)) =\mathbb E(K(X_1,.)\ell(Y_1)).
 \end{displaymath}
 \item For every $K,K'\in\mathcal K_n$,
 \begin{displaymath}
 \mathbb E(\langle K(X_1,.),K'(X_2,.)\ell(Y_2)\rangle_{2}^{2})
 \leqslant\mathfrak m_{\mathcal K,\ell}\overline s_{K',\ell}
 \end{displaymath}
 with
 \begin{displaymath}
 \overline s_{K',\ell} :=\mathbb E(\|K'(X_1,.)\ell(Y_1)\|_{2}^{2}).
 \end{displaymath}
 \item For every $K\in\mathcal K_n$ and $\psi\in\mathbb L^2(\mathbb R^d)$,
 \begin{displaymath}
 \mathbb E(\langle K(X_1,.),\psi\rangle_{2}^{2})\leqslant\mathfrak m_{\mathcal K,\ell}\|\psi\|_{2}^{2}.
 \end{displaymath}
\end{enumerate}
\end{assumption}
\noindent
The elements of $\mathcal K_n$ are called kernels. Let us provide two natural examples of kernels sets.
%


%
\begin{proposition}\label{example_kernels_set_1}
Consider
\begin{displaymath}
\mathcal K_k(h_{\min}) :=
\left\{
(x',x)\mapsto
\prod_{q = 1}^{d}\frac{1}{h_q}k\left(\frac{x_q' - x_q}{h_q}\right)
\textrm{ $;$ }
h_1,\dots, h_d\in\mathcal H(h_{\min})\right\},
\end{displaymath}
where $k$ is a symmetric kernel (in the usual sense), $h_{\min}\in [n^{-1/d},1]$ and $\mathcal H(h_{\min})$ is a finite subset of $[h_{\min},1]$. The kernels set $\mathcal K_k(h_{\min})$ fulfills Assumption \ref{assumption_K} and, for any $K\in\mathcal K_k(h_{\min})$ (i.e. defined by (\ref{PR_kernel}) with $h_1,\dots,h_d\in\mathcal H(h_{\min})$),
\begin{displaymath}
\overline s_{K,\ell} =\|k\|_{2}^{2d}\mathbb E(\ell(Y_1)^2)
\prod_{q = 1}^{d}\frac{1}{h_q}.
\end{displaymath}
\end{proposition}
%


%
\begin{proposition}\label{example_kernels_set_2}
Consider
\begin{displaymath}
\mathcal K_{\mathcal B_1,\dots,\mathcal B_n}(m_{\max}) :=
\left\{
(x',x)\mapsto
\prod_{q = 1}^{d}
\sum_{j = 1}^{m_q}\varphi_{j}^{m_q}(x_q)\varphi_{j}^{m_q}(x_q')
\textrm{ $;$ }
m_1,\dots,m_d\in\{1,\dots,m_{\max}\}\right\},
\end{displaymath}
where $m_{\max}^{d}\in\{1,\dots,n\}$ and, for every $m\in\{1,\dots,n\}$, $\mathcal B_m =\{\varphi_{1}^{m},\dots,\varphi_{m}^{m}\}$ is an orthonormal family of $\mathbb L^2(\mathbb R)$ such that
\begin{displaymath}
\sup_{x'\in\mathbb R}\sum_{j = 1}^{m}\varphi_{j}^{m}(x')^2
\leqslant\mathfrak m_{\mathcal B}m
\end{displaymath}
with $\mathfrak m_{\mathcal B} > 0$ not depending on $m$ and $n$, and such that one of the two following conditions is satisfied:
\begin{equation}\label{projection_condition_1}
\mathcal B_m\subset\mathcal B_{m + 1}
\textrm{ $;$ }
\forall m\in\{1,\dots,n - 1\}
\end{equation}
or
\begin{equation}\label{projection_condition_2}
\overline{\mathfrak m}_{\mathcal B} :=
\sup\{
|\mathbb E(K(X_1,x))|\textrm{ $;$ }
K\in\mathcal K_{\mathcal B_1,\dots,\mathcal B_n}(m_{\max})
\textrm{ and }x\in\mathbb R^d\}
\textrm{ is finite and doesn't depend on $n$.}
\end{equation}
The kernels set $\mathcal K_{\mathcal B_1,\dots,\mathcal B_n}(m_{\max})$ fulfills Assumption \ref{assumption_K} and, for any $K\in\mathcal K_{\mathcal B_1,\dots,\mathcal B_n}(m_{\max})$ (i.e. defined by (\ref{projection_kernel}) with $m_1,\dots,m_n\in\{1,\dots,m_{\max}\}$),
\begin{displaymath}
\overline s_{K,\ell}\leqslant
\mathfrak m_{\mathcal B}^{d}\mathbb E(\ell(Y_1)^2)
\prod_{q = 1}^{d}m_q.
\end{displaymath}
\end{proposition}
%


%
\begin{remark}\label{weighted_projection_kernels_set}
For the sake of simplicity, the present paper focuses on $\mathcal K_{\mathcal B_1,\dots,\mathcal B_n}(m_{\max})$, but Proposition \ref{example_kernels_set_2} is still true for the weighted projection kernels set
\begin{displaymath}
\mathcal K_{\mathcal B_1,\dots,\mathcal B_n}(w_1,\dots,w_n;m_{\max}) :=
\left\{
(x',x)\mapsto
\prod_{q = 1}^{d}
\sum_{j = 1}^{m_q}w_j\varphi_{j}^{m_q}(x_q)\varphi_{j}^{m_q}(x_q')
\textrm{ $;$ }
m_1,\dots,m_d\in\{1,\dots,m_{\max}\}\right\},
\end{displaymath}
where $w_1,\dots,w_n\in [0,1]$.
\end{remark}
%


%
\begin{remark}\label{remark_projection_conditions}
Note that Condition (\ref{projection_condition_1}) is close, but more restrictive than Condition (19) of Lerasle et al. \cite{LMRB16}, Proposition 3.2, which is that the spaces $\normalfont{\textrm{span}}(\mathcal B_m)$, $m\in\mathbb N$ are nested. See Massart \cite{MASSART07}, Subsection 7.5.2 for examples of nested spaces. Our Condition (\ref{projection_condition_1}) is fulfilled by the trigonometric basis, Hermite's basis or Laguerre's basis.\\
Note also that in the same proposition of Lerasle et al. \cite{LMRB16}, Condition (20) coincides with our Condition (\ref{projection_condition_2}). The regular histograms basis satisfies Condition (\ref{projection_condition_2}). Indeed, by taking $\varphi_{j}^{m} =\psi_{j}^{m} :=\sqrt m\mathbf 1_{[(j - 1)/m,j/m[}$ for every $m\in\{1,\dots,n\}$ and $j\in\{1,\dots,m\}$,
\begin{eqnarray*}
 \left|\mathbb E\left[
 \prod_{q = 1}^{d}
 \sum_{j = 1}^{m_q}\psi_{j}^{m_q}(X_{1,q})\psi_{j}^{m_q}(x_q)\right]\right| & = &
 \sum_{j_1 = 1}^{m_1}\cdots
 \sum_{j_d = 1}^{m_d}
 \left(\prod_{q = 1}^{d}m_q\mathbf 1_{[(j_q - 1)/m_q,j_q/m_q[}(x_q)\right)\\
 & &
 \quad\quad\times
 \int_{(j_1 - 1)/m_1}^{j_1/m_1}\cdots\int_{(j_d - 1)/m_d}^{j_d/m_d}f(x_1',\dots,x_d')dx_1'\cdots dx_d'\\
 & \leqslant &
 \|f\|_{\infty}
 \prod_{q = 1}^{d}\sum_{j = 1}^{m_q}\mathbf 1_{[(j - 1)/m_q,j/m_q[}(x_q)\leqslant\|f\|_{\infty}
\end{eqnarray*}
for every $m_1,\dots,m_d\in\{1,\dots,n\}$ and $x\in\mathbb R^d$.
\end{remark}
\noindent
The following proposition shows that Legendre's basis also fulfills Condition (\ref{projection_condition_2}).
%


%
\begin{proposition}\label{Legendre_basis}
For every $m\in\{1,\dots,n\}$ and $j\in\{1,\dots,m\}$, let $\xi_{j}^{m}$ be the function defined on $[-1,1]$ by
\begin{displaymath}
\xi_{j}^{m}(x) :=
\sqrt{\frac{2j + 1}{2}}Q_j(x)
\textrm{ $;$ }
\forall x\in [-1,1],
\end{displaymath}
where
\begin{displaymath}
Q_j : x\in [-1,1]\longmapsto
\frac{1}{2^jj!}\cdot\frac{d^j}{dx^j}(x^2 - 1)^j
\end{displaymath}
is the $j$-th Legendre's polynomial. If $f\in C^{2d}([0,1]^d)$ and $\mathcal B_m =\{\xi_{1}^{m},\dots,\xi_{m}^{m}\}$ for every $m\in\{1,\cdots,n\}$, then $\mathcal K_{\mathcal B_1,\dots,\mathcal B_m}(m_{\max})$ fulfills Condition (\ref{projection_condition_2}).
\end{proposition}
\noindent
The following proposition provides a suitable control of the variance of $\widehat s_{K,\ell}(n;.)$.
%


%
\begin{proposition}\label{variance_bound_main_estimator}
Under Assumption \ref{assumption_K}.(1,2,3), if $s\in\mathbb L^2(\mathbb R^d)$ and if there exists $\alpha > 0$ such that $\mathbb E(\exp(\alpha|\ell(Y_1)|)) <\infty$, then there exists a deterministic constant $\mathfrak c_{\ref{variance_bound_main_estimator}} > 0$, not depending on $n$, such that for every $\theta\in ]0,1[$,
\begin{displaymath}
\mathbb E\left(\sup_{K\in\mathcal K_n}\left\{
\left|\|\widehat s_{K,\ell}(n;.) - s_{K,\ell}\|_{2}^{2} -\frac{\overline s_{K,\ell}}{n}\right|
-\frac{\theta}{n}\overline s_{K,\ell}\right\}\right)
\leqslant
\mathfrak c_{\ref{variance_bound_main_estimator}}\frac{\log(n)^5}{\theta n}.
\end{displaymath}
\end{proposition}
\noindent
Finally, let us state the main result of this section.
%


%
\begin{theorem}\label{risk_bound_main_estimator}
Under Assumption \ref{assumption_K}, if $s\in\mathbb L^2(\mathbb R^d)$ and if there exists $\alpha > 0$ such that $\mathbb E(\exp(\alpha|\ell(Y_1)|)) <\infty$, then there exist deterministic constants $\mathfrak c_{\ref{risk_bound_main_estimator}},\overline{\mathfrak c}_{\ref{risk_bound_main_estimator}} > 0$, not depending on $n$, such that for every $\theta\in ]0,1[$,
\begin{displaymath}
\mathbb E\left(\sup_{K\in\mathcal K_n}\left\{
\|\widehat s_{K,\ell}(n;.) - s\|_{2}^{2} - (1 +\theta)\left(\|s_{K,\ell} - s\|_{2}^{2} +\frac{\overline s_{K,\ell}}{n}\right)\right\}\right)
\leqslant
\mathfrak c_{\ref{risk_bound_main_estimator}}\frac{\log(n)^5}{\theta n}
\end{displaymath}
and
\begin{displaymath}
\mathbb E\left(\sup_{K\in\mathcal K_n}\left\{
\|s_{K,\ell} - s\|_{2}^{2} +\frac{\overline s_{K,\ell}}{n} -
\frac{1}{1 -\theta}\|\widehat s_{K,\ell}(n;.) - s\|_{2}^{2}\right\}\right)
\leqslant
\overline{\mathfrak c}_{\ref{risk_bound_main_estimator}}\frac{\log(n)^5}{\theta(1 -\theta) n}.
\end{displaymath}
\end{theorem}
%


%
\begin{remark}\label{remark_risk_bound_main_estimator_1}
Note that the first inequality in Theorem \ref{risk_bound_main_estimator} gives a risk bound on the estimator $\widehat s_{K,\ell}(n;.)$:
\begin{displaymath}
\mathbb E(\|\widehat s_{K,\ell}(n;.) - s\|_{2}^{2})
\leqslant
(1 +\theta)\left(\|s_{K,\ell} - s\|_{2}^{2} +\frac{\overline s_{K,\ell}}{n}\right) +
\mathfrak c_{\ref{risk_bound_main_estimator}}
\frac{\log(n)^5}{\theta n}
\end{displaymath}
for every $\theta\in ]0,1[$. The second inequality is useful in order to establish a risk bound on the adaptive estimator defined in the next section (see Theorem \ref{risk_bound_adaptive_estimator}).
\end{remark}
%


%
\begin{remark}\label{remark_risk_bound_main_estimator_2}
In Proposition \ref{variance_bound_main_estimator} and Theorem \ref{risk_bound_main_estimator}, the exponential moment condition may appear too strong. Nevertheless, this is de facto satisfied when
\begin{equation}\label{remark_risk_bound_main_estimator_2_1}
\ell(Y_1),\dots,\ell(Y_n)\textrm{ have a compactly supported distribution.}
\end{equation}
This last condition is satisfied in the density estimation framework because $\ell = 1$, but even in the nonparametric regression framework, where $\ell$ is not bounded, when $Y_1,\dots,Y_n$ have a compactly supported distribution. Moreover, note that under Condition (\ref{remark_risk_bound_main_estimator_2_1}), the risk bounds of Theorem \ref{risk_bound_main_estimator} can be stated in deviation, without additional steps in the proof. Precisely, under Assumption \ref{assumption_K} and Condition (\ref{remark_risk_bound_main_estimator_2_1}), if $s\in\mathbb L^2(\mathbb R^d)$, then there exists a deterministic constant $\mathfrak c_L > 0$, depending on $L =\sup_{z\in\normalfont{\textrm{supp}}(\mathbb P_{Y_1})}|\ell(z)|$ but not on $n$, such that for every $\vartheta\in]0,1[$ and $\lambda > 0$,
\begin{displaymath}
\sup_{K\in\mathcal K_n}\left|
\|s_{K,\ell} - s\|_{2}^{2} +\frac{\overline s_{K,\ell}}{n} -
\frac{1}{1 -\vartheta}\|\widehat s_{K,\ell}(n;.) - s\|_{2}^{2}\right|\leqslant\frac{\mathfrak c_L}{\vartheta n}(1+\lambda)^3
\end{displaymath}
with probability larger than $1 - 9.4|\mathcal K_n|e^{-\lambda}$.\\
When Condition (\ref{remark_risk_bound_main_estimator_2_1}) doesn't hold true, one can replace the exponential moment condition of Proposition \ref{variance_bound_main_estimator} and Theorem \ref{risk_bound_main_estimator} by a $q$-th order moment condition on $\ell(Y_1)$ ($q\in\mathbb N^*$), but with a damaging effect on the rate of convergence of $\widehat s_{K,\ell}(n;.)$. For instance, at Remark \ref{remark_bound_trace_term}, it is established that under a $(12 - 4\varepsilon)/\beta$-th moment condition ($\varepsilon\in ]0,1[$ and $0 <\beta < \varepsilon/2$), the rate of convergence is of order $O(1/n^{1 -\varepsilon})$ (instead of $1/n$) in Lemma \ref{bound_trace_term}. This holds true for the three technical lemmas of Subsection \ref{subsection_preliminary_results}, and then for Proposition \ref{variance_bound_main_estimator} and Theorem \ref{risk_bound_main_estimator}.
\end{remark}
%


%
\section{Kernel selection}\label{section_kernel_selection}
This section deals with a risk bound on the adaptive estimator $\widehat s_{\widehat K,\ell}(n;.)$, where
\begin{displaymath}
\widehat K\in\arg\min_{K\in\mathcal K_n}
\{\|\widehat s_{K,\ell}(n;\cdot) -\widehat s_{K_0,\ell}(n;\cdot)\|_{2}^{2} +\textrm{pen}_{\ell}(K)\},
\end{displaymath}
$K_0$ is an overfitting proposal for $K$ in the sense that
\begin{displaymath}
K_0\in\arg\max_{K\in\mathcal K_n}\left\{
\sup_{x\in\mathbb R^d}|K(x,x)|\right\},
\end{displaymath}
and
\begin{equation}\label{penalty_proposal}
\textrm{pen}_{\ell}(K) :=
\frac{2}{n^2}
\sum_{i = 1}^{n}\langle K(.,X_i),K_0(.,X_i)\rangle_2\ell(Y_i)^2
\textrm{ $;$ }
\forall K\in\mathcal K_n.
\end{equation}
\textbf{Example.} On the one hand, for any $K\in\mathcal K_k(h_{\min})$ (i.e. defined by (\ref{PR_kernel}) with $h_1,\dots,h_d\in\mathcal H(h_{\min})$),
\begin{displaymath}
\sup_{x\in\mathbb R^d}|K(x,x)| =
|k(0)|^d\prod_{q = 1}^{d}\frac{1}{h_q}.
\end{displaymath}
Then, for $\mathcal K_n =\mathcal K_k(h_{\min})$,
\begin{displaymath}
K_0(x',x) =
\frac{1}{h_{\min}^{d}}\prod_{q = 1}^{d}k\left(\frac{x_q' - x_q}{h_{\min}}\right)
\textrm{ $;$ }
\forall x,x'\in\mathbb R^d.
\end{displaymath}
On the other hand, for any $K\in\mathcal K_{\mathcal B_1,\dots,\mathcal B_n}(m_{\max})$ (i.e. defined by (\ref{projection_kernel}) with $m_1,\dots,m_n\in\{1,\dots,m_{\max}\}$),
\begin{displaymath}
\sup_{x\in\mathbb R^d}|K(x,x)| =
\sup_{x\in\mathbb R^d}
\prod_{q = 1}^{d}\sum_{j = 1}^{m_q}\varphi_{j}^{m_q}(x_q)^2.
\end{displaymath}
Then, for $\mathcal K_n =\mathcal K_{\mathcal B_1,\dots,\mathcal B_n}(m_{\max})$, at least for the usual bases mentioned at Remark \ref{remark_projection_conditions},
\begin{displaymath}
K_0(x',x) =
\prod_{q = 1}^{d}
\sum_{j = 1}^{m_{\max}}
\varphi_{j}^{m_{\max}}(x_q)\varphi_{j}^{m_{\max}}(x_q')
\textrm{ $;$ }
\forall x,x'\in\mathbb R^d.
\end{displaymath}
In the sequel, in addition to Assumption \ref{assumption_K}, the kernels set $\mathcal K_n$ fulfills the following assumption.
%


%
\begin{assumption}\label{additional_assumption_K}
There exists a deterministic constant $\overline{\mathfrak m}_{\mathcal K,\ell} > 0$, not depending on $n$, such that
\begin{displaymath}
\mathbb E\left(\sup_{K,K'\in\mathcal K_n}
\langle K(X_1,.),s_{K',\ell}\rangle_{2}^{2}\right)
\leqslant\overline{\mathfrak m}_{\mathcal K,\ell}.
\end{displaymath}
\end{assumption}
\noindent
The following theorem provides an oracle inequality for the adaptive estimator $\widehat s_{\widehat K,\ell}(n;.)$.
%


%
\begin{theorem}\label{risk_bound_adaptive_estimator}
Under Assumptions \ref{assumption_K} and \ref{additional_assumption_K}, if $s\in\mathbb L^2(\mathbb R^d)$ and if there exists $\alpha > 0$ such that $\mathbb E(\exp(\alpha|\ell(Y_1)|)) <\infty$, then there exists a deterministic constant $\mathfrak c_{\ref{risk_bound_adaptive_estimator}} > 0$, not depending on $n$, such that for every $\vartheta\in ]0,1[$,
\begin{displaymath}
\mathbb E(\|\widehat s_{\widehat K,\ell}(n;.) - s\|_{2}^{2})
\leqslant
(1 +\vartheta)\min_{K\in\mathcal K_n}
\mathbb E(\|\widehat s_{K,\ell}(n;.) - s\|_{2}^{2}) +
\frac{\mathfrak c_{\ref{risk_bound_adaptive_estimator}}}{\vartheta}
\left(\|s_{K_0,\ell} - s\|_{2}^{2} +\frac{\log(n)^5}{n}\right).
\end{displaymath}
\end{theorem}
%


%
\begin{remark}\label{rate_convergence_PCO_estimator}
As mentioned in Comte and Marie \cite{CM}, p. 6, when $\mathcal K_n =\mathcal K_k(h_{\min})$, if $s$ belongs to a Nikol'skii ball and $h_{\min} = 1/n$, then Theorem \ref{risk_bound_adaptive_estimator} says that the PCO estimator has a performance of same order than $O_n :=\min_{K\in\mathcal K_n}\mathbb E(\|\widehat s_{K,\ell}(n;.) - s\|_{2}^{2})$ up to a factor $1 +\vartheta$. When $\mathcal K_n =\mathcal K_{\mathcal B_1,\dots,\mathcal B_n}(m_{\max})$, it depends on the bases $\mathcal B_1,\dots,\mathcal B_n$. For instance, with the same ideas than in Comte and Marie \cite{CM}, thanks to DeVore and Lorentz \cite{DL93}, Theorem 2.3 p. 205, if $s$ belongs to a Sobolev space and $m_{\max} = n$, then our Theorem \ref{risk_bound_adaptive_estimator} also says that the PCO estimator has a performance of same order than $O_n$.
\end{remark}
\noindent
\textbf{Notation.} For any $B\in\mathcal B(\mathbb R^d)$, $\|.\|_{2,f,B}$ is the norm on $\mathbb L^2(B,f(x)\lambda_d(dx))$ defined by
\begin{displaymath}
\|\varphi\|_{2,f,B} :=
\left(\int_B\varphi(x)^2f(x)\lambda_d(dx)\right)^{1/2}
\textrm{$;$ }
\forall\varphi\in\mathbb L^2(B,f(x)\lambda_d(dx)).
\end{displaymath}
The following corollary provides an oracle inequality for $\widehat s_{\widehat K,\ell}(n;.)/\widehat s_{\widehat L,1}(n;.)$, where $\widehat K$ and $\widehat L$ are both selected via the PCO method.
%


%
\begin{corollary}\label{risk_bound_adaptive_quotient_estimator}
Let $(\beta_j)_{j\in\mathbb N}$ be a decreasing sequence of elements of $]0,\infty[$ such that $\lim_{\infty}\beta_j = 0$ and, for every $j\in\mathbb N$, consider
\begin{displaymath}
B_j :=
\{x\in\mathbb R^d : f(x)\geqslant\beta_j\}.
\end{displaymath}
Under Assumptions \ref{assumption_K} and \ref{additional_assumption_K} for $\ell$ and $1$, if $s,f\in\mathbb L^2(\mathbb R^d)$ and if there exists $\alpha > 0$ such that $\mathbb E(\exp(\alpha|\ell(Y_1)|)) <\infty$, then there exists a deterministic constant $\mathfrak c_{\ref{risk_bound_adaptive_estimator}} > 0$, not depending on $n$, such that for every $\vartheta\in ]0,1[$,
\begin{eqnarray*}
 \mathbb E\left[\left\|\frac{\widehat s_{\widehat K,\ell}(n;.)}{\widehat s_{\widehat L,1}(n;.)} -\frac{s}{f}\right\|_{2,f,B_n}^{2}\right]
 & \leqslant &
 \frac{\mathfrak c_{\ref{risk_bound_adaptive_estimator}}}{\beta_{n}^{2}}\left[
 (1 +\vartheta)\min_{(K,L)\in\mathcal K_{n}^{2}}\{
 \mathbb E(\|\widehat s_{K,\ell}(n;.) - s\|_{2}^{2}) +\mathbb E(\|\widehat s_{L,1}(n;.) - f\|_{2}^{2})\}
 \right.\\
 & &
 \left.
 +\frac{1}{\vartheta}
 \left(\|s_{K_0,\ell} - s\|_{2}^{2} +\|s_{K_0,1} - f\|_{2}^{2} +\frac{\log(n)^5}{n}\right)\right]
\end{eqnarray*}
where
\begin{displaymath}
\widehat K\in\arg\min_{K\in\mathcal K_n}
\{\|\widehat s_{K,\ell}(n;\cdot) -\widehat s_{K_0,\ell}(n;\cdot)\|_{2}^{2} +
\normalfont{\textrm{pen}}_{\ell}(K)\}
\end{displaymath}
and
\begin{displaymath}
\widehat L\in\arg\min_{L\in\mathcal K_n}
\{\|\widehat s_{L,1}(n;\cdot) -\widehat s_{K_0,1}(n;\cdot)\|_{2}^{2} +
\normalfont{\textrm{pen}}_1(L)\}.
\end{displaymath}
\end{corollary}
\noindent
The proof of Corollary \ref{risk_bound_adaptive_quotient_estimator} is the same than the proof of Comte and Marie \cite{CM}, Corollary 4.3.
\\
\\
Finally, let us discuss about Assumption \ref{additional_assumption_K}. This assumption is difficult to check in practice, then let us provide a sufficient condition.
%


%
\begin{assumption}\label{additional_assumption_K_sufficient}
The function $s$ is bounded and
\begin{displaymath}
\mathfrak m_{\mathcal K} :=
\sup\{\|K(x',.)\|_{1}^{2}
\textrm{ $;$ }K\in\mathcal K_n\textrm{ and }x'\in\mathbb R^d\}
\end{displaymath}
doesn't depend on $n$.
\end{assumption}
\noindent
Under Assumption \ref{additional_assumption_K_sufficient}, $\mathcal K_n$ fulfills Assumption \ref{additional_assumption_K}. Indeed,
\begin{small}
\begin{eqnarray*}
 \mathbb E\left(\sup_{K,K'\in\mathcal K_n}
 \langle K(X_1,.),s_{K',\ell}\rangle_{2}^{2}\right)
 & \leqslant &
 \left(\sup_{K'\in\mathcal K_n}\|s_{K',\ell}\|_{\infty}^{2}\right)\mathbb E\left(\sup_{K\in\mathcal K_n}\|K(X_1,.)\|_{1}^{2}\right)\\
 & \leqslant &
 \mathfrak m_{\mathcal K}\sup\left\{
 \left(\int_{-\infty}^{\infty}|K'(x',x)s(x)|dx\right)^2
 \textrm{ $;$ }
 K'\in\mathcal K_n\textrm{ and }x'\in\mathbb R\right\}
 \leqslant
 \mathfrak m_{\mathcal K}^{2}\|s\|_{\infty}^{2}.
\end{eqnarray*}
\end{small}
\newline
Note that in the nonparametric regression framework (see Model (\ref{nonparametric_regression})), to assume $s$ bounded means that $bf$ is bounded. For instance, this condition is fulfilled by the linear regression models with Gaussian inputs.
\\
Let us provide two examples of kernels sets fulfilling Assumption \ref{additional_assumption_K_sufficient}, the sufficient condition for Assumption \ref{additional_assumption_K}:
\begin{itemize}
 \item Consider $K\in\mathcal K_k(h_{\min})$. Then, there exist $h_1,\dots,h_d\in\mathcal H(h_{\min})$ such that
 \begin{displaymath}
 K(x',x) =
 \prod_{q = 1}^{d}
 \frac{1}{h_q}k\left(\frac{x_q' - x_q}{h_q}\right)
 \textrm{ $;$ }
 \forall x,x'\in\mathbb R^d.
 \end{displaymath}
 Clearly, $\|K(x',.)\|_1 =\|k\|_{1}^{d}$ for every $x'\in\mathbb R^d$. So, for $\mathcal K_n =\mathcal K_k(h_{\min})$, $\mathfrak m_{\mathcal K}\leqslant\|k\|_{1}^{2d}$.
 \item For $\mathcal K_n =\mathcal K_{\mathcal B_1,\dots,\mathcal B_n}(m_{\max})$, the condition on $\mathfrak m_{\mathcal K}$ seems harder to check in general. Let us show that it is satisfied for the regular histograms basis defined in Section \ref{section_risk_bound}. For every $m_1,\dots,m_d\in\{1,\dots,n\}$,
 \begin{displaymath}
 \left\|\prod_{q = 1}^{d}\sum_{j = 1}^{m_q}\psi_{j}^{m_q}(x_q')\psi_{j}^{m_q}(.)\right\|_1
 \leqslant
 \prod_{q = 1}^{d}\left(
 m_q\sum_{j = 1}^{m_q}\mathbf 1_{[(j - 1)/m_q,j/m_q[}(x_q')\int_{(j - 1)/m_q}^{j/m_q}dx\right)
 \leqslant 1.
 \end{displaymath}
\end{itemize}
Now, let us show that even if it doesn't fulfill Assumption \ref{additional_assumption_K_sufficient}, the trigonometric basis fulfills Assumption \ref{additional_assumption_K}.
%


%
\begin{proposition}\label{example_kernels_set+}
Consider $\chi_1 :=\mathbf 1_{[0,1]}$ and, for every $j\in\mathbb N^*$, the functions $\chi_{2j}$ and $\chi_{2j + 1}$ defined on $\mathbb R$ by
\begin{displaymath}
\chi_{2j}(x) :=\sqrt 2\cos(2\pi jx)\mathbf 1_{[0,1]}(x)
\textrm{ and }
\chi_{2j + 1}(x) :=\sqrt 2\sin(2\pi j x)\mathbf 1_{[0,1]}(x)
\textrm{ $;$ }
\forall x\in\mathbb R.
\end{displaymath}
If $s\in C^2(\mathbb R^d)$ and $\mathcal B_m =\{\chi_1,\dots,\chi_m\}$ for every $m\in\{1,\dots,n\}$, then $\mathcal K_{\mathcal B_1,\dots,\mathcal B_n}(m_{\max})$ fulfills Assumption \ref{additional_assumption_K}.
\end{proposition}
\appendix
%


%
\section{Details on kernels sets: proofs of Propositions \ref{example_kernels_set_1}, \ref{example_kernels_set_2}, \ref{Legendre_basis} and \ref{example_kernels_set+}}
%


%
\subsection{Proof of Proposition \ref{example_kernels_set_1}}
Consider $K,K'\in\mathcal K_k(h_{\min})$. Then, there exist $h,h'\in\mathcal H(h_{\min})^d$ such that
\begin{displaymath}
K(x',x) =
k_h(x' - x)
\textrm{ and }
K'(x',x) =
k_{h'}(x' - x)
\end{displaymath}
for every $x,x'\in\mathbb R^d$, where
\begin{displaymath}
k_h(x) :=
\prod_{q = 1}^{d}
\frac{1}{h_q}k\left(\frac{x_q}{h_q}\right)
\textrm{ $;$ }
\forall x\in\mathbb R^d.
\end{displaymath}
\begin{enumerate}
 \item For every $x'\in\mathbb R^d$, since $nh_{\min}^{d}\geqslant 1$,
 \begin{eqnarray}
  \label{example_kernels_set_1_1}
  \|K(x',.)\|_{2}^{2} & = &
  \left(\prod_{q = 1}^{d}
  \frac{1}{h_{q}^{2}}\right)\left[
  \int_{\mathbb R^d}\prod_{q = 1}^{d}
  k\left(\frac{x_q' - x_q}{h_q}\right)^2\lambda_d(dx)\right]
  =\|k\|_{2}^{2d}\prod_{q = 1}^{d}\frac{1}{h_q}\\
  & \leqslant &
  \|k\|_{2}^{2d}\frac{1}{h_{\min}^{d}}
  \leqslant\|k\|_{2}^{2d}n.
  \nonumber
 \end{eqnarray}
 \item Since $s_{K,\ell} = K\ast s$ and by Young's inequality, $\|s_{K,\ell}\|_{2}^{2}\leqslant\|k\|_{1}^{2d}\|s\|_{2}^{2}$.
 \item On the one hand, thanks to Equality (\ref{example_kernels_set_1_1}),
 \begin{displaymath}
 \overline s_{K',\ell} =
 \mathbb E(\|K'(X_1,.)\ell(Y_1)\|_{2}^{2}) =
 \|k\|_{2}^{2d}\mathbb E(\ell(Y_1)^2)\prod_{q = 1}^{d}
 \frac{1}{h_q'}.
 \end{displaymath}
 On the other hand, for every $x,x'\in\mathbb R^d$,
 \begin{displaymath}
 \langle K(x,.),K'(x',.)\rangle_2 =
 \int_{\mathbb R^d}k_h(x - x'')k_{h'}(x' - x'')\lambda_d(dx'') =
 (k_h\ast k_{h'})(x - x').
 \end{displaymath}
 Then,
 \begin{eqnarray*}
  \mathbb E(\langle K(X_1,.),K'(X_2,.)\ell(Y_2)\rangle_{2}^{2})
  & = &
  \mathbb E((k_h\ast k_{h'})(X_1 - X_2)^2\ell(Y_2)^2)\\
  & = &
  \int_{\mathbb R^{d + 1}}\left[\ell(y)^2\int_{\mathbb R^d}(k_h\ast k_{h'})(x' - x)^2f(x')\lambda_d(dx')\right]
  \mathbb P_{(X_2,Y_2)}(dx,dy)\\
  & \leqslant &
  \|f\|_{\infty}\|k_h\ast k_{h'}\|_{2}^{2}
  \mathbb E(\ell(Y_2)^2)
  \leqslant
  \|f\|_{\infty}\|k\|_{1}^{2d}\overline s_{K',\ell}.
 \end{eqnarray*}
 \item For every $\psi\in\mathbb L^2(\mathbb R^d)$,
 \begin{eqnarray*}
  \mathbb E(\langle K(X_1,.),\psi\rangle_{2}^{2})
  & = &
  \mathbb E((k_h\ast\psi)(X_1)^2)\\
  & \leqslant &
  \|f\|_{\infty}
  \|k_h\ast\psi\|_{2}^{2}
  \leqslant
  \|f\|_{\infty}
  \|k\|_{1}^{2d}\|\psi\|_{2}^{2}.
 \end{eqnarray*}
\end{enumerate}
%


%
\subsection{Proof of Proposition \ref{example_kernels_set_2}}
Consider $K,K'\in\mathcal K_{\mathcal B_1,\dots,\mathcal B_n}(m_{\max})$. Then, there exist $m,m'\in\{1,\dots,m_{\max}\}^d$ such that
\begin{displaymath}
K(x',x) =
\prod_{q = 1}^{d}
\sum_{j = 1}^{m_q}\varphi_{j}^{m_q}(x_q)\varphi_{j}^{m_q}(x_q')
\textrm{ and }
K'(x',x) =
\prod_{q = 1}^{d}
\sum_{j = 1}^{m_q'}\varphi_{j}^{m_q'}(x_q)\varphi_{j}^{m_q'}(x_q')
\end{displaymath}
for every $x,x'\in\mathbb R^d$.
\begin{enumerate}
 \item For every $x'\in\mathbb R^d$, since $m_{\max}^{d}\leqslant n$,
 \begin{eqnarray}
  \label{example_kernels_set_2_1}
  \|K(x',.)\|_{2}^{2} & = &
  \prod_{q = 1}^{d}
  \sum_{j,j' = 1}^{m_q}\varphi_{j'}^{m_q}(x_q')\varphi_{j}^{m_q}(x_q')
  \int_{-\infty}^{\infty}\varphi_{j'}^{m_q}(x)\varphi_{j}^{m_q}(x)dx
  =\prod_{q = 1}^{d}
  \sum_{j = 1}^{m_q}
  \varphi_{j}^{m_q}(x_q')^2\\
  & \leqslant &
  \mathfrak m_{\mathcal B}^{d}
  \prod_{q = 1}^{d}m_q
  \leqslant
  \mathfrak m_{\mathcal B}^{d}n.
  \nonumber
 \end{eqnarray}
 \item Since
 \begin{displaymath}
 s_{K,\ell}(.) =
 \sum_{j_1 = 1}^{m_1}\cdots\sum_{j_d = 1}^{m_d}
 \langle s,\varphi_{j_1}^{m_1}\otimes\cdots\otimes\varphi_{j_d}^{m_d}\rangle_2
 (\varphi_{j_1}^{m_1}\otimes\cdots\otimes\varphi_{j_d}^{m_d})(.),
 \end{displaymath}
 by Pythagoras theorem, $\|s_{K,\ell}\|_{2}^{2}\leqslant\|s\|_{2}^{2}$.
 \item First of all, thanks to Equality (\ref{example_kernels_set_2_1}),
 \begin{displaymath}
 \overline s_{K',\ell}
 = \mathbb E\left[\ell(Y_1)^2
 \prod_{q = 1}^{d}
 \sum_{j = 1}^{m_q'}\varphi_{j}^{m_q'}(X_{1,q})^2\right]
 \leqslant\mathfrak m_{\mathcal B}^{d}\mathbb E(\ell(Y_1)^2)
 \prod_{q = 1}^{d}m_q'.
 \end{displaymath}
 On the one hand, under Condition (\ref{projection_condition_1}) on $\mathcal B_1,\dots,\mathcal B_n$, for any $j\in\{1,\dots,m\}$, $\varphi_{j}^{m}$ doesn't depend on $m$, so it can be denoted by $\varphi_j$, and then
 \begin{eqnarray*}
  \mathbb E(\langle K(X_1,.),K'(X_2,.)\ell(Y_2)\rangle_{2}^{2})
  & = &
  \int_{\mathbb R^d}\mathbb E\left[
  \left(\prod_{q = 1}^{d}
  \sum_{j = 1}^{m_q\wedge m_q'}\varphi_j(x_q')\varphi_j(X_{2,q})\right)^2\ell(Y_2)^2\right]f(x')\lambda_d(dx')\\
  & \leqslant &
  \|f\|_{\infty}\mathbb E\left[\ell(Y_2)^2
  \prod_{q = 1}^{d}
  \sum_{j,j' = 1}^{m_q\wedge m_q'}\varphi_{j'}(X_{2,q})\varphi_j(X_{2,q})
  \int_{-\infty}^{\infty}\varphi_{j'}(x')\varphi_j(x')dx'\right]\\
  & &
  \quad\quad\leqslant
  \|f\|_{\infty}\overline s_{K',\ell}.
 \end{eqnarray*}
 On the other hand, under Condition (\ref{projection_condition_2}) on $\mathcal B_1,\dots,\mathcal B_n$, since $X_1$ and $(X_2,Y_2)$ are independent, and since $K(x,x)\geqslant 0$ for every $x\in\mathbb R^d$,
 \begin{eqnarray*}
  \mathbb E(\langle K(X_1,.),K'(X_2,.)\ell(Y_2)\rangle_{2}^{2})
  & \leqslant &
  \mathbb E(\|K(X_1,.)\|_{2}^{2}\|K'(X_2,.)\|_{2}^{2}\ell(Y_2)^2)\\
  & = & \mathbb E(K(X_1,X_1))
  \mathbb E(\|K'(X_2,.)\|_{2}^{2}\ell(Y_2)^2)
  \leqslant
  \overline{\mathfrak m}_{\mathcal B}\overline s_{K',\ell}.
 \end{eqnarray*}
 \item For every $\psi\in\mathbb L^2(\mathbb R^d)$,
 \begin{eqnarray*}
  \mathbb E(\langle K(X_1,.),\psi\rangle_{2}^{2})
  & = &
  \mathbb E\left[\left|\sum_{j_1 = 1}^{m_1}\cdots\sum_{j_d = 1}^{m_d}
  \langle\psi,\varphi_{j_1}^{m_1}\otimes\cdots\otimes\varphi_{j_d}^{m_d}\rangle_2
  (\varphi_{j_1}^{m_1}\otimes\cdots\otimes\varphi_{j_d}^{m_d})(X_1)\right|^2\right]\\
  & \leqslant &
  \|f\|_{\infty}\left\|\sum_{j_1 = 1}^{m_1}\cdots\sum_{j_d = 1}^{m_d}
  \langle\psi,\varphi_{j_1}^{m_1}\otimes\cdots\otimes\varphi_{j_d}^{m_d}\rangle_2
  (\varphi_{j_1}^{m_1}\otimes\cdots\otimes\varphi_{j_d}^{m_d})(.)\right\|_{2}^{2}
  \leqslant\|f\|_{\infty}\|\psi\|_{2}^{2}.
 \end{eqnarray*}
\end{enumerate}
%


%
\subsection{Proof of Proposition \ref{Legendre_basis}}
For the sake of readability, assume that $d = 1$. Consider $m\in\{1,\dots,m_{\max}\}$. Since each Legendre's polynomial is uniformly bounded by $1$,
\begin{displaymath}
\left|\mathbb E\left[\sum_{j = 1}^{m}\xi_{j}^{m}(X_1)\xi_{j}^{m}(x')\right]\right|
\leqslant\sum_{j = 1}^{m}\frac{2j + 1}{2}\left|\int_{-1}^{1}Q_j(x)f(x)dx\right|.
\end{displaymath}	
Moreover, since $Q_j$ is a solution to Legendre's differential equation for any $j\in\{1,\dots,m\}$, thanks to the integration by parts formula,
\begin{eqnarray*}
 \int_{-1}^{1}Q_j(x)f(x)dx
 & = & -\frac{1}{j(j + 1)}\int_{-1}^{1}\frac{d}{dx}[(1 - x^2)Q_j'(x)]f(x)dx\\
 & = & -\frac{1}{j(j + 1)}[(1 - x^2)Q_j'(x)f(x)]_{-1}^{1}
 +\frac{1}{j(j + 1)}\int_{-1}^{1}(1 - x^2)Q_j'(x)f'(x)dx\\	
 & = & -\frac{1}{j(j + 1)}\int_{-1}^{1}Q_j(x)\frac{d}{dx}[(1 - x^2)f'(x)]dx.
\end{eqnarray*}
Then,
\begin{displaymath}
\left|\int_{-1}^{1}Q_j(x)f(x)dx\right|
\leqslant \frac{2\mathfrak c_1}{j(j + 1)}\|Q_j\|_2
=\frac{2\sqrt 2\mathfrak c_1}{j(j + 1)(2j + 1)^{1/2}}
\end{displaymath}
with $\mathfrak c_1 =\max\{2\|f'\|_{\infty},\|f''\|_{\infty}\}$. So,
\begin{displaymath}
\left|\mathbb E\left[\sum_{j = 1}^{m}\xi_{j}^{m}(X_1)\xi_{j}^{m}(x')\right]\right|
\leqslant
2\mathfrak c_1\sum_{j = 1}^{m}
\frac{1}{j^{3/2}}
\leqslant 2\mathfrak c_1\zeta\left(\frac{3}{2}\right)
\end{displaymath}
where $\zeta$ is Riemann's zeta function. Thus, Legendre's basis satisfies Condition (\ref{projection_condition_2}).
%


%
\subsection{Proof of Proposition \ref{example_kernels_set+}}
The proof of Proposition \ref{example_kernels_set+} relies on the following technical lemma.
%


%
\begin{lemma}\label{lemma_example_kernels_set+}
For every $x\in [0,2\pi]$ and $p,q\in\mathbb N^*$ such that $q > p$,
\begin{displaymath}
\left|\sum_{j = p + 1}^{q}\frac{\sin(jx)}{j}\right|
\leqslant
\frac{2}{(1 + p)\sin(x/2)}.
\end{displaymath}
\end{lemma}
\noindent
See Subsubsection \ref{proof_lemma_example_kernels_set+} for a proof.
\\
\\
For the sake of readability, assume that $d = 1$. Consider $K,K'\in\mathcal K_{\mathcal B_1,\dots,\mathcal B_n}(m_{\max})$. Then, there exist $m,m'\in\{1,\dots,m_{\max}\}$ such that
\begin{displaymath}
K(x',x) =
\sum_{j = 1}^{m}\chi_j(x)\chi_j(x')
\textrm{ and }
K'(x',x) =
\sum_{j = 1}^{m'}\chi_j(x)\chi_j(x')
\textrm{ $;$ }
\forall x,x'\in\mathbb R.
\end{displaymath}
First, there exist $\mathfrak m_1(m,m')\in\{0,\dots,n\}$ and $\mathfrak c_1 > 0$, not depending on $n$, $K$ and $K'$, such that for any $x'\in [0,1]$,
\begin{eqnarray*}
 |\langle K(x',.),s_{K',\ell}\rangle_2| & = &
 \left|\sum_{j = 1}^{m\wedge m'}
 \mathbb E(\ell(Y_1)\chi_j(X_1))\chi_j(x')\right|\\
 & \leqslant &
 \mathfrak c_1 + 2\left|
 \sum_{j = 1}^{\mathfrak m_1(m,m')}
 \mathbb E(\ell(Y_1)(\cos(2\pi jX_1)\cos(2\pi jx') +
 \sin(2\pi jX_1)\sin(2\pi jx'))\mathbf 1_{[0,1]}(X_1))\right|\\
 & = &
 \mathfrak c_1 + 2\left|
 \sum_{j = 1}^{\mathfrak m_1(m,m')}
 \mathbb E(\ell(Y_1)\cos(2\pi j(X_1 - x'))\mathbf 1_{[0,1]}(X_1))\right|. 
\end{eqnarray*}
Moreover, for any $j\in\{2,\dots,\mathfrak m_1(m,m')\}$,
\begin{eqnarray*}
 \mathbb E(\ell(Y_1)\cos(2\pi j(X_1 - x'))\mathbf 1_{[0,1]}(X_1)) & = &
 \int_{0}^{1}\cos(2\pi j(x - x'))s(x)dx\\
 & = &
 \frac{1}{j}\left[\frac{\sin(2\pi j(x - x'))}{2\pi}s(x)\right]_{0}^{1}\\
 & &
 +\frac{1}{j^2}\left[\frac{\cos(2\pi j(x - x'))}{4\pi^2}s'(x)\right]_{0}^{1}
 -\frac{1}{j^2}\int_{0}^{1}\frac{\cos(2\pi j(x - x'))}{4\pi^2}s''(x)dx\\
 & = &
 \frac{s(0) - s(1)}{2\pi}\cdot\frac{\alpha_j(x')}{j} +\frac{\beta_j(x')}{j^2}
\end{eqnarray*}
where $\alpha_j(x') :=\sin(2\pi jx')$ and
\begin{displaymath}
\beta_j(x') :=
\frac{1}{4\pi^2}\left(
(s'(1) - s'(0))\cos(2\pi jx') -\int_{0}^{1}\cos(2\pi j(x - x'))s''(x)dx\right).
\end{displaymath}
Then, there exists a deterministic constant $\mathfrak c_2 > 0$, not depending on $n$, $K$, $K'$ and $x'$, such that
\begin{equation}\label{example_kernels_set+_1}
\langle K(x',.),s_{K',\ell}\rangle_{2}^{2}
\leqslant
\mathfrak c_2\left[1 +
\left(\sum_{j = 1}^{\mathfrak m_1(m,m')}\frac{\alpha_j(x')}{j}\right)^2 +
\left(\sum_{j = 1}^{\mathfrak m_1(m,m')}\frac{\beta_j(x')}{j^2}\right)^2\right].
\end{equation}
Let us show that each term of the right-hand side of Inequality (\ref{example_kernels_set+_1}) is uniformly bounded in $x'$, $m$ and $m'$. On the one hand,
\begin{displaymath}
\left|\sum_{j = 1}^{\mathfrak m_1(m,m')}\frac{\beta_j(x')}{j^2}\right|
\leqslant
\max_{j\in\{1,\dots,n\}}
\|\beta_j\|_{\infty}
\sum_{j = 1}^{n}\frac{1}{j^2}
\leqslant
\frac{1}{24}
(2\|s'\|_{\infty} +\|s''\|_{\infty}).
\end{displaymath}
On the other hand, for every $x\in ]0,\pi[$ such that $[\pi/x] + 1\leqslant\mathfrak m_1(m,m')$ (without loss of generality), by Lemma \ref{lemma_example_kernels_set+},
\begin{eqnarray}
 \left|\sum_{j = 1}^{\mathfrak m_1(m,m')}\frac{\sin(jx)}{j}\right|
 & \leqslant &
 \left|\sum_{j = 1}^{[\pi/x]}\frac{\sin(jx)}{j}\right| +
 \left|\sum_{j = [\pi/x] + 1}^{\mathfrak m_1(m,m')}\frac{\sin(jx)}{j}\right|
 \nonumber\\
 \label{example_kernels_set+_2}
 & \leqslant & x\left[\frac{\pi}{x}\right] +\frac{2}{(1 + [\pi/x])\sin(x/2)}
 \leqslant\pi + 2.
\end{eqnarray}
Since $x\mapsto\sin(x)$ is continuous, odd and $2\pi$-periodic, Inequality (\ref{example_kernels_set+_2}) holds true for every $x\in\mathbb R$. So,
\begin{displaymath}
\left|\sum_{j = 1}^{\mathfrak m_1(m,m')}\frac{\alpha_j(x')}{j}\right|
\leqslant
\pi + 2.
\end{displaymath}
Therefore,
\begin{displaymath}
\mathbb E\left[\sup_{K,K'\in\mathcal K_{\mathcal B_1,\dots,\mathcal B_n}(m_{\max})}
\langle K(X_1,.),s_{K',\ell}\rangle_{2}^{2}\right]
\leqslant
\mathfrak c_2\left(1 +(\pi + 2)^2 +
\frac{1}{24^2}
(2\|s'\|_{\infty} +\|s''\|_{\infty})^2\right).
\end{displaymath}
%


%
\subsubsection{Proof of Lemma \ref{lemma_example_kernels_set+}}\label{proof_lemma_example_kernels_set+}
For any $x\in [0,2\pi]$ and $q\in\mathbb N^*$, consider
\begin{displaymath}
f_q(x) :=
\sum_{j = 1}^{q}\frac{\sin(jx)}{j}
\textrm{, }
g_q(x) :=
\sum_{j = 1}^{q}\left(\frac{1}{j} -\frac{1}{j + 1}\right)h_j(x)\textrm{ and }
h_q(x) :=
\sum_{j = 1}^{q}\sin(jx).
\end{displaymath}
On the one hand,
\begin{displaymath}
g_q(x) = h_1(x) -\frac{1}{q + 1}h_q(x) +\sum_{j = 2}^{q}\frac{1}{j}(h_j(x) - h_{j - 1}(x)).
\end{displaymath}
Then,
\begin{displaymath}
f_q(x) = g_q(x) +\frac{1}{q + 1}h_q(x).
\end{displaymath}
On the other hand,
\begin{eqnarray*}
 h_q(x) & = &
 \textrm{Im}\left(\sum_{j = 1}^{q}e^{\mathbf ijx}\right) =
 \textrm{Im}\left[e^{\mathbf i(q + 1)x/2}\frac{\sin(qx/2)}{\sin(x/2)}\right]\\
 & = &
 \frac{\sin((q + 1)x/2)\sin(qx/2)}{\sin(x/2)}
 =\frac{\cos(x/2) -\cos((q + 1/2)x)}{2\sin(x/2)}.
\end{eqnarray*}
Then,
\begin{displaymath}
\sin\left(\frac{x}{2}\right)|h_q(x)|
\leqslant
1
\end{displaymath}
and, for any $p\in\mathbb N^*$ such that $q > p$,
\begin{displaymath}
\sin\left(\frac{x}{2}\right)|g_q(x) - g_p(x)|
\leqslant\frac{1}{p + 1} -\frac{1}{q + 1}.
\end{displaymath}
Therefore,
\begin{eqnarray*}
 \sin\left(\frac{x}{2}\right)
 |f_q(x) - f_p(x)|
 & \leqslant &
 \sin\left(\frac{x}{2}\right)|g_q(x) - g_p(x)| +
 \sin\left(\frac{x}{2}\right)\frac{|h_q(x)|}{q + 1} +
 \sin\left(\frac{x}{2}\right)\frac{|h_p(x)|}{p + 1}\\
 & \leqslant &
 \frac{2}{p + 1}.
\end{eqnarray*}
In conclusion,
\begin{displaymath}
\left|\sum_{j = p + 1}^{q}\frac{\sin(jx)}{k}\right|
\leqslant
\frac{2}{(1 + p)\sin(x/2)}.
\end{displaymath}
%


%
\section{Proofs of risk bounds}
%


%
\subsection{Preliminary results}\label{subsection_preliminary_results}
This subsection provides three lemmas used several times in the sequel.
%


%
\begin{lemma}\label{bound_U_statistics}
Consider
\begin{equation}\label{U_statistic_definition}
U_{K,K',\ell}(n) :=\sum_{i\not= j}
\langle K(X_i,.)\ell(Y_i) - s_{K,\ell},
K'(X_j,.)\ell(Y_j) - s_{K',\ell}\rangle_2
\textrm{ $;$ }
\forall K,K'\in\mathcal K_n.
\end{equation}
Under Assumption \ref{assumption_K}.(1,2,3), if $s\in\mathbb L^2(\mathbb R^d)$ and if there exists $\alpha > 0$ such that $\mathbb E(\exp(\alpha|\ell(Y_1)|)) <\infty$, then there exists a deterministic constant $\mathfrak c_{\ref{bound_U_statistics}} > 0$, not depending on $n$, such that for every $\theta\in ]0,1[$,
\begin{displaymath}
\mathbb E\left(\sup_{K,K'\in\mathcal K_n}
\left\{\frac{|U_{K,K',\ell}(n)|}{n^2}
-\frac{\theta}{n}\overline s_{K',\ell}
\right\}\right)
\leqslant
\mathfrak c_{\ref{bound_U_statistics}}
\frac{\log(n)^5}{\theta n}.
\end{displaymath}
\end{lemma}
%


%
\begin{lemma}\label{bound_trace_term}
Consider
\begin{displaymath}
V_{K,\ell}(n) :=\frac{1}{n}\sum_{i = 1}^{n}\|K(X_i,.)\ell(Y_i) - s_{K,\ell}\|_{2}^{2}
\textrm{ $;$ }
\forall K\in\mathcal K_n.
\end{displaymath}
Under Assumption \ref{assumption_K}.(1,2), if $s\in\mathbb L^2(\mathbb R^d)$ and if there exists $\alpha > 0$ such that $\mathbb E(\exp(\alpha|\ell(Y_1)|)) <\infty$, then there exists a deterministic constant $\mathfrak c_{\ref{bound_trace_term}} > 0$, not depending on $n$, such that for every $\theta\in ]0,1[$,
\begin{displaymath}
\mathbb E\left(\sup_{K\in\mathcal K_n}\left\{
\frac{1}{n}|V_{K,\ell}(n) -\overline s_{K,\ell}| -\frac{\theta}{n}\overline s_{K,\ell}\right\}\right)
\leqslant\mathfrak c_{\ref{bound_trace_term}}\frac{\log(n)^3}{\theta n}.
\end{displaymath}
\end{lemma}
%


%
\begin{lemma}\label{bound_crossed_term}
Consider
\begin{equation}\label{W_statistic_definition}
W_{K,K',\ell}(n) :=
\langle\widehat s_{K,\ell}(n;.) - s_{K,\ell},s_{K',\ell} - s\rangle_2
\textrm{ $;$ }
\forall K,K'\in\mathcal K_n.
\end{equation}
Under Assumption \ref{assumption_K}.(1,2,4), if $s\in\mathbb L^2(\mathbb R^d)$ and if there exists $\alpha > 0$ such that $\mathbb E(\exp(\alpha|\ell(Y_1)|)) <\infty$, then there exists a deterministic constant $\mathfrak c_{\ref{bound_crossed_term}} > 0$, not depending on $n$, such that for every $\theta\in ]0,1[$,
\begin{displaymath}
\mathbb E\left(\sup_{K,K'\in\mathcal K_n}\{
|W_{K,K',\ell}(n)| -\theta\|s_{K',\ell} - s\|_{2}^{2}
\}\right)
\leqslant
\mathfrak c_{\ref{bound_crossed_term}}\frac{\log(n)^4}{\theta n}.
\end{displaymath}
\end{lemma}
%


%
\subsubsection{Proof of Lemma \ref{bound_U_statistics}}
The proof of Lemma \ref{bound_U_statistics} relies on the following concentration inequality for U-statistics, proved in dimension $1$ in Houdr\'e and Reynaud-Bouret \cite{HRB03} first, and then extended to the infinite-dimensional framework by Gin\'e and Nickl in \cite{GN15}.
%


%
\begin{lemma}\label{U_statistic_concentration_inequality}
Let $\xi_1,\dots,\xi_n$ be i.i.d. random variables on a Polish space $\Xi$ equipped with its Borel $\sigma$-algebra. Let $f_{i,j}$, $1\leqslant i\not= j\leqslant n$, be some bounded and symmetric measurable maps from $\Xi^2$ into $\mathbb R$ such that, for every $i\not= j$,
\begin{displaymath}
f_{i,j} = f_{j,i}
\textrm{ and }
\mathbb E(f_{i,j}(z,\xi_1)) = 0\textrm{ $dz$-a.e.}
\end{displaymath}
Consider the totally degenerate second order U-statistic
\begin{displaymath}
U_n :=
\sum_{i\not= j}f_{i,j}(\xi_i,\xi_j).
\end{displaymath}
There exists a universal constant $\mathfrak m > 0$ such that for every $\lambda > 0$,
\begin{displaymath}
\mathbb P(U_n\leqslant\mathfrak m(\mathfrak c_n\lambda^{1/2} +\mathfrak d_n\lambda +\mathfrak b_n\lambda^{3/2} +\mathfrak a_n\lambda^2))\geqslant 1 - 2.7e^{-\lambda}
\end{displaymath}
where
\begin{eqnarray*}
 \mathfrak a_n & = &
 \sup_{i,j = 1,\dots,n}\left\{\sup_{z,z'\in\Xi}|f_{i,j}(z,z')|\right\},\\
 \mathfrak b_{n}^{2} & = &
 \max\left\{\sup_{i,z}\sum_{j = 1}^{i - 1}\mathbb E(f_{i,j}(z,\xi_j)^2)\textrm{ $;$ }
 \sup_{j,z'}\sum_{i = j + 1}^{n}\mathbb E(f_{i,j}(\xi_i,z')^2)\right\},\\
 \mathfrak c_{n}^{2} & = &
 \sum_{i\not= j}\mathbb E(f_{i,j}(\xi_i,\xi_j)^2)\textrm{ and}\\
 \mathfrak d_n & = &
 \sup_{(a,b)\in\mathcal A}\mathbb E\left[\sum_{i < j}f_{i,j}(\xi_i,\xi_j)a_i(\xi_i)b_j(\xi_j)\right]
\end{eqnarray*}
with
\begin{displaymath}
\mathcal A =
\left\{(a,b) :
\mathbb E\left(\sum_{i = 1}^{n - 1}a_i(\xi_i)^2\right)\leqslant 1
\textrm{ and }
\mathbb E\left(\sum_{j = 2}^{n}b_j(\xi_j)^2\right)\leqslant 1\right\}.
\end{displaymath}
\end{lemma}
\noindent
See Gin\'e and Nickl \cite{GN15}, Theorem 3.4.8 for a proof.
\\
\\
Consider $\mathfrak m(n) := 8\log(n)/\alpha$. For any $K,K'\in\mathcal K_n$,
\begin{displaymath}
U_{K,K',\ell}(n)
= U_{K,K',\ell}^{1}(n) + U_{K,K',\ell}^{2}(n) + U_{K,K',\ell}^{3}(n) + U_{K,K',\ell}^{4}(n)
\end{displaymath}
where
\begin{displaymath}
U_{K,K',\ell}^{l}(n) :=
\sum_{i \not= j}
g_{K,K',\ell}^{l}(n;X_i,Y_i,X_j,Y_j)
\textrm{ $;$ }
l = 1,2,3,4
\end{displaymath}
with, for every $(x',y),(x'',y')\in E =\mathbb R^d\times\mathbb R$,
\begin{eqnarray*}
 g_{K,K',\ell}^{1}(n;x',y,x'',y') & := &
 \langle K(x',.)\ell(y)\mathbf 1_{|\ell(y)|\leqslant\mathfrak m(n)} - s_{K,\ell}^{+}(n;.),
 K'(x'',.)\ell(y')\mathbf 1_{|\ell(y)|\leqslant\mathfrak m(n)} - s_{K',\ell}^{+}(n;.)\rangle_2,\\
 g_{K,K',\ell}^{2}(n;x',y,x'',y') & := &
 \langle K(x',.)\ell(y)\mathbf 1_{|\ell(y)| >\mathfrak m(n)} - s_{K,\ell}^{-}(n;.),
 K'(x'',.)\ell(y')\mathbf 1_{|\ell(y)|\leqslant\mathfrak m(n)} - s_{K',\ell}^{+}(n;.)\rangle_2,\\
 g_{K,K',\ell}^{3}(n;x',y,x'',y') & := &
 \langle K(x',.)\ell(y)\mathbf 1_{|\ell(y)|\leqslant\mathfrak m(n)} - s_{K,\ell}^{+}(n;.),
 K'(x'',.)\ell(y')\mathbf 1_{|\ell(y)| >\mathfrak m(n)} - s_{K',\ell}^{-}(n;.)\rangle_2,\\
 g_{K,K',\ell}^{4}(n;x',y,x'',y') & := &
 \langle K(x',.)\ell(y)\mathbf 1_{|\ell(y)| >\mathfrak m(n)} - s_{K,\ell}^{-}(n;.),
 K'(x'',.)\ell(y')\mathbf 1_{|\ell(y)| >\mathfrak m(n)} - s_{K',\ell}^{-}(n;.)\rangle_2
\end{eqnarray*}
and, for every $k\in\mathcal K_n$,
\begin{displaymath}
s_{k,\ell}^{+}(n;.) :=
\mathbb E(k(X_1,.)\ell(Y_1)\mathbf 1_{|\ell(Y_1)|\leqslant\mathfrak m(n)})
\textrm{ and }
s_{k,\ell}^{-}(n;.) :=
\mathbb E(k(X_1,.)\ell(Y_1)\mathbf 1_{|\ell(Y_1)| >\mathfrak m(n)}).
\end{displaymath}
On the one hand, since $\mathbb E(g_{K,K',\ell}^{1}(n;x',y,X_1,Y_1)) = 0$ for every $(x',y)\in E$, by Lemma \ref{U_statistic_concentration_inequality}, there exists a universal constant $\mathfrak m\geqslant 1$ such that for any $\lambda > 0$, with probability larger than $1 - 5.4e^{-\lambda}$,
\begin{displaymath}
\frac{|U_{K,K',\ell}^{1}(n)|}{n^2}
\leqslant\frac{\mathfrak m}{n^2}(\mathfrak c_{K,K',\ell}(n)\lambda^{1/2} +\mathfrak d_{K,K',\ell}(n)\lambda +\mathfrak b_{K,K',\ell}(n)\lambda^{3/2} +\mathfrak a_{K,K',\ell}(n)\lambda^2)
\end{displaymath}
where the constants $\mathfrak a_{K,K',\ell}(n)$, $\mathfrak b_{K,K',\ell}(n)$, $\mathfrak c_{K,K',\ell}(n)$ and $\mathfrak d_{K,K',\ell}(n)$ are defined and controlled later. First, note that
\begin{eqnarray}
 U_{K,K',\ell}^{1}(n) & = &
 \sum_{i\not= j}
 (\varphi_{K,K',\ell}(n;X_i,Y_i,X_j,Y_j)
 \nonumber\\
 \label{bound_U_statistics_1}
 & &
 \quad\quad\quad
 -\psi_{K,K',\ell}(n;X_i,Y_i)
 -\psi_{K',K,\ell}(n;X_j,Y_j)
 +\mathbb E(\varphi_{K,K',\ell}(n;X_i,Y_i,X_j,Y_j))),
\end{eqnarray}
where
\begin{displaymath}
\varphi_{K,K',\ell}(n;x',y,x'',y'') :=
\langle K(x',.)\ell(y)\mathbf 1_{|\ell(y)|\leqslant\mathfrak m(n)},
K'(x'',.)\ell(y')\mathbf 1_{|\ell(y')|\leqslant\mathfrak m(n)}\rangle_2
\end{displaymath}
and
\begin{eqnarray*}
\psi_{k,k',\ell}(n;x',y) :=
\langle k(x',.)\ell(y)\mathbf 1_{|\ell(y)|\leqslant\mathfrak m(n)},
s_{k',\ell}^{+}(n;.)\rangle_2 =
\mathbb E(\varphi_{k,k',\ell}(n;x',y,X_1,Y_1))
\end{eqnarray*}
for every $k,k'\in\mathcal K_n$ and $(x',y),(x'',y')\in E$. Let us now control $\mathfrak a_{K,K',\ell}(n)$, $\mathfrak b_{K,K',\ell}(n)$, $\mathfrak c_{K,K',\ell}(n)$ and $\mathfrak d_{K,K',\ell}(n)$:
\begin{itemize}
 \item\textbf{The constant $\mathfrak a_{K,K',\ell}(n)$.} Consider
 \begin{displaymath}
 \mathfrak a_{K,K',\ell}(n) :=
 \sup_{(x',y),(x'',y')\in E}
 |g_{K,K',\ell}^{1}(n;x',y,x'',y')|.
 \end{displaymath}
 By (\ref{bound_U_statistics_1}), Cauchy-Schwarz's inequality and Assumption \ref{assumption_K}.(1),
 \begin{eqnarray*}
  \mathfrak a_{K,K',\ell}(n) & \leqslant &
  4\sup_{(x',y),(x'',y')\in E}
  |\langle K(x',.)\ell(y)\mathbf 1_{|\ell(y)|\leqslant\mathfrak m(n)},
  K'(x'',.)\ell(y')\mathbf 1_{|\ell(y')|\leqslant\mathfrak m(n)}\rangle_2|\\
  & \leqslant &
  4\mathfrak m(n)^2
  \left(\sup_{x'\in\mathbb R^d}\|K(x',.)\|_2\right)
  \left(\sup_{x''\in\mathbb R^d}\|K'(x'',.)\|_2\right)
  \leqslant 4\mathfrak m_{\mathcal K,\ell}\mathfrak m(n)^2n.
 \end{eqnarray*}
 So,
 \begin{displaymath}
 \frac{1}{n^2}\mathfrak a_{K,K',\ell}(n)\lambda^2
 \leqslant
 \frac{4}{n}\mathfrak m_{\mathcal K,\ell}\mathfrak m(n)^2\lambda^2.
 \end{displaymath}
 \item\textbf{The constant $\mathfrak b_{K,K',\ell}(n)$.} Consider
 \begin{displaymath}
 \mathfrak b_{K,K',\ell}(n)^2 :=
 n\sup_{(x',y)\in E}
 \mathbb E(g_{K,K',\ell}^{1}(n;x',y,X_1,Y_1)^2).
 \end{displaymath}
 By (\ref{bound_U_statistics_1}), Jensen's inequality, Cauchy-Schwarz's inequality and Assumption \ref{assumption_K}.(1),
 \begin{eqnarray*}
  \mathfrak b_{K,K',\ell}(n)^2 & \leqslant &
  16n\sup_{(x',y)\in E}
  \mathbb E(\langle K(x',.)\ell(y)\mathbf 1_{|\ell(y)|\leqslant\mathfrak m(n)},
  K'(X_1,.)\ell(Y_1)\mathbf 1_{|\ell(Y_1)|\leqslant\mathfrak m(n)}\rangle_{2}^{2})\\
  & \leqslant &
  16n\mathfrak m(n)^2\sup_{x'\in\mathbb R^d}
  \|K(x',.)\|_{2}^{2}
  \mathbb E(
  \|K'(X_1,.)\ell(Y_1)\mathbf 1_{|\ell(Y_1)|\leqslant\mathfrak m(n)}\|_{2}^{2})
  \leqslant
  16\mathfrak m_{\mathcal K,\ell}n^2\mathfrak m(n)^2\overline s_{K',\ell}.
 \end{eqnarray*}
 So, for any $\theta\in]0,1[$,
 \begin{eqnarray*}
  \frac{1}{n^2}\mathfrak b_{K,K',\ell}(n)\lambda^{3/2}
  & \leqslant &
  2\left(\frac{3\mathfrak m}{\theta}\right)^{1/2}\frac{2}{n^{1/2}}\mathfrak m_{\mathcal K,\ell}^{1/2}\mathfrak m(n)\lambda^{3/2}
  \times\left(\frac{\theta}{3\mathfrak m}\right)^{1/2}\frac{1}{n^{1/2}}\overline s_{K',\ell}^{1/2}\\
  & \leqslant &
  \frac{\theta}{3\mathfrak mn}\overline s_{K',\ell} +
  \frac{12\mathfrak m\lambda^3}{\theta n}\mathfrak m_{\mathcal K,\ell}\mathfrak m(n)^2.
 \end{eqnarray*}
 \item\textbf{The constant $\mathfrak c_{K,K',\ell}(n)$.} Consider
 \begin{displaymath}
 \mathfrak c_{K,K',\ell}(n)^2 :=
 n^2\mathbb E(g_{K,K',\ell}^{1}(n;X_1,Y_1,X_2,Y_2)^2).
 \end{displaymath}
 By (\ref{bound_U_statistics_1}), Jensen's inequality and Assumption \ref{assumption_K}.(3),
 \begin{eqnarray*}
  \mathfrak c_{K,K',\ell}(n)^2
  & \leqslant &
  16n^2\mathbb E(\langle K(X_1,.)\ell(Y_1)\mathbf 1_{|\ell(Y_1)|\leqslant\mathfrak m(n)},
  K'(X_2,.)\ell(Y_2)\mathbf 1_{|\ell(Y_2)|\leqslant\mathfrak m(n)}\rangle_{2}^{2})\\
  & \leqslant &
  16n^2\mathfrak m(n)^2
  \mathbb E(\langle K(X_1,.),K'(X_2,.)\ell (Y_2)\rangle_{2}^{2})
  \leqslant 16\mathfrak m_{\mathcal K,\ell}n^2\mathfrak m(n)^2\overline s_{K',\ell}.
 \end{eqnarray*}
 So,
 \begin{displaymath}
 \frac{1}{n^2}\mathfrak c_{K,K',\ell}(n)\lambda^{1/2}
 \leqslant
 \frac{\theta}{3\mathfrak mn}\overline s_{K',\ell} +
 \frac{12\mathfrak m\lambda}{\theta n}\mathfrak m_{\mathcal K,\ell}\mathfrak m(n)^2.
 \end{displaymath}
 \item\textbf{The constant $\mathfrak d_{K,K',\ell}(n)$.} Consider
 \begin{displaymath}
 \mathfrak d_{K,K',\ell}(n) :=
 \sup_{(a,b)\in\mathcal A}
 \mathbb E\left[\sum_{i < j}a_i(X_i,Y_i)b_j(X_j,Y_j)g_{K,K',\ell}^{1}(n;X_i,Y_i,X_j,Y_j)\right],
 \end{displaymath}
 where
 \begin{displaymath}
 \mathcal A :=
 \left\{(a,b) :
 \sum_{i = 1}^{n - 1}\mathbb E(a_i(X_i,Y_i)^2)\leqslant 1
 \textrm{ and }
 \sum_{j = 2}^{n}\mathbb E(b_j(X_j,Y_j)^2)\leqslant 1\right\}.
 \end{displaymath}
 By (\ref{bound_U_statistics_1}), Jensen's inequality, Cauchy-Schwarz's inequality and Assumption \ref{assumption_K}.(3),
 \begin{eqnarray*}
  \mathfrak d_{K,K',\ell}(n) & \leqslant &
  4\sup_{(a,b)\in\mathcal A}\mathbb E\left[
  \sum_{i = 1}^{n - 1}
  \sum_{j = i + 1}^{n}
  |a_i(X_i,Y_i)b_j(X_j,Y_j)
  \varphi_{K,K',\ell}(n;X_i,Y_i,X_j,Y_j)|\right]\\
  & \leqslant &
  4n\mathfrak m(n)
  \mathbb E(\langle K(X_1,.),K'(X_2,.)\ell (Y_2)\rangle_{2}^{2})^{1/2}
  \leqslant
  4\mathfrak m_{\mathcal K,\ell}^{1/2}n\mathfrak m(n)\overline s_{K',\ell}^{1/2}.
 \end{eqnarray*}
 So,
 \begin{displaymath}
 \frac{1}{n^2}\mathfrak d_{K,K',\ell}(n)\lambda
 \leqslant
 \frac{\theta}{3\mathfrak mn}\overline s_{K',\ell} +
 \frac{12\mathfrak m\lambda^2}{\theta n}\mathfrak m_{\mathcal K,\ell}\mathfrak m(n)^2.
 \end{displaymath}
\end{itemize}
Then, since $\mathfrak m\geqslant 1$ and $\lambda > 0$, with probability larger than $1 - 5.4e^{-\lambda}$,
\begin{displaymath}
\frac{|U_{K,K',\ell}^{1}(n)|}{n^2}
\leqslant
\frac{\theta}{n}\overline s_{K',\ell} +
\frac{40\mathfrak m^2}{\theta n}\mathfrak m_{\mathcal K,\ell}\mathfrak m(n)^2(1 +\lambda)^3.
\end{displaymath}
So, with probability larger than $1 - 5.4|\mathcal K_n|^2e^{-\lambda}$,
\begin{displaymath}
S_{\mathcal K,\ell}(n,\theta) :=
\sup_{K,K'\in\mathcal K_n}
\left\{\frac{|U_{K,K',\ell}^{1}(n)|}{n^2}
-\frac{\theta}{n}\overline s_{K',\ell}
\right\}
\leqslant
\frac{40\mathfrak m^2}{\theta n}\mathfrak m_{\mathcal K,\ell}\mathfrak m(n)^2(1 +\lambda)^3.
\end{displaymath}
For every $t\in\mathbb R_+$, consider
\begin{displaymath}
\lambda_{\mathcal K,\ell}(n,\theta,t) :=
-1 +\left(\frac{t}{\mathfrak m_{\mathcal K,\ell}(n,\theta)}\right)^{1/3}
\textrm{ with }
\mathfrak m_{\mathcal K,\ell}(n,\theta) =
\frac{40\mathfrak m^2}{\theta n}\mathfrak m_{\mathcal K,\ell}\mathfrak m(n)^2.
\end{displaymath}
Then, for any $T > 0$,
\begin{eqnarray*}
 \mathbb E(S_{\mathcal K,\ell}(n,\theta))
 & \leqslant &
 T +\int_{T}^{\infty}\mathbb P(S_{\mathcal K,\ell}(n,\theta)\geqslant
 (1 +\lambda_{\mathcal K,\ell}(n,\theta,t))^3\mathfrak m_{\mathcal K,\ell}(n,\theta))dt\\
 & \leqslant &
 T + 5.4|\mathcal K_n|^2\int_{T}^{\infty}\exp(-\lambda_{\mathcal K,\ell}(n,\theta,t))dt\\
 & = &
 T + 5.4|\mathcal K_n|^2\int_{T}^{\infty}
 \exp\left(-\frac{t^{1/3}}{2\mathfrak m_{\mathcal K,\ell}(n,\theta)^{1/3}}\right)
 \exp\left(1 -\frac{t^{1/3}}{2\mathfrak m_{\mathcal K,\ell}(n,\theta)^{1/3}}\right)dt\\
 & \leqslant &
 T + 5.4\mathfrak c_1|\mathcal K_n|^2\mathfrak m_{\mathcal K,\ell}(n,\theta)
 \exp\left(-\frac{T^{1/3}}{2\mathfrak m_{\mathcal K,\ell}(n,\theta)^{1/3}}\right)
 \textrm{ with }
 \mathfrak c_1 =\int_{0}^{\infty}e^{1 - r^{1/3}/2}dr.
\end{eqnarray*}
Moreover,
\begin{displaymath}
\mathfrak m_{\mathcal K,\ell}(n,\theta)
\leqslant\mathfrak c_2\frac{\log(n)^2}{\theta n}
\textrm{ with }
\mathfrak c_2 =
\frac{40\cdot 8^2\mathfrak m^2}{\alpha^2}\mathfrak m_{\mathcal K,\ell}.
\end{displaymath}
So, by taking
\begin{displaymath}
T = 2^4\mathfrak c_2\frac{\log(n)^5}{\theta n},
\end{displaymath}
and since $|\mathcal K_n|\leqslant n$,
\begin{displaymath}
\mathbb E(S_{\mathcal K,\ell}(n,\theta))
\leqslant
2^4\mathfrak c_2\frac{\log(n)^5}{\theta n} +
5.4\mathfrak c_1\mathfrak m_{\mathcal K,\ell}(n,\theta)
\frac{|\mathcal K_n|^2}{n^2}
\leqslant
(2^4 + 5.4\mathfrak c_1)
\mathfrak c_2\frac{\log(n)^5}{\theta n}.
\end{displaymath}
On the other hand, by Assumption \ref{assumption_K}.(1), Cauchy-Schwarz's inequality and Markov's inequality,
\begin{small}
\begin{eqnarray*}
 \mathbb E\left(\sup_{K,K'\in\mathcal K_n}
 |g_{K,K',\ell}^{2}(n ; X_1,Y_1,X_2,Y_2)|\right) & \leqslant &
 4\mathfrak m(n)\sum_{K,K'\in\mathcal K_n}
 \mathbb E(|\ell(Y_1)|\mathbf 1_{|\ell(Y_1)| >\mathfrak m(n)}
 |\langle K(X_1,.),K'(X_2,.)\rangle_2|)\\
 & \leqslant &
 4\mathfrak m(n)\mathfrak m_{\mathcal K,\ell}n|\mathcal K_n|^2
 \mathbb E(\ell(Y_1)^2)^{1/2}
 \mathbb P(|\ell(Y_1)| >\mathfrak m(n))^{1/2}
 \leqslant
 \mathfrak c_3
 \frac{\log(n)}{n}
\end{eqnarray*}
\end{small}
\newline
with
\begin{displaymath}
\mathfrak c_3 =
\frac{32}{\alpha}\mathfrak m_{\mathcal K,\ell}
\mathbb E(\ell(Y_1)^2)^{1/2}
\mathbb E(\exp(\alpha|\ell(Y_1)|))^{1/2}.
\end{displaymath}
So,
\begin{displaymath}
\mathbb E\left(\sup_{K,K'\in\mathcal K_n}
\frac{|U_{K,K',\ell}^{2}(n)|}{n^2}\right)
\leqslant
\mathfrak c_3
\frac{\log(n)}{n}
\end{displaymath}
and, symmetrically,
\begin{displaymath}
\mathbb E\left(\sup_{K,K'\in\mathcal K_n}
\frac{|U_{K,K',\ell}^{3}(n)|}{n^2}\right)
\leqslant
\mathfrak c_3
\frac{\log(n)}{n}.
\end{displaymath}
By Assumption \ref{assumption_K}.(1), Cauchy-Schwarz's inequality and Markov's inequality,
\begin{small}
\begin{eqnarray*}
 \mathbb E\left(\sup_{K,K'\in\mathcal K_n}
 |g_{K,K',\ell}^{4}(n ; X_1,Y_1,X_2,Y_2)|\right) & \leqslant &
 4\sum_{K,K'\in\mathcal K_n}
 \mathbb E(|\ell(Y_1)\ell(Y_2)|\mathbf 1_{|\ell(Y_1)|,|\ell(Y_2)| >\mathfrak m(n)}
 |\langle K(X_1,.),K'(X_2,.)\rangle_2|)\\
 & \leqslant &
 4\mathfrak m_{\mathcal K,\ell}n|\mathcal K_n|^2
 \mathbb E(\ell(Y_1)^2)
 \mathbb P(|\ell(Y_1)| >\mathfrak m(n))
 \leqslant
 \frac{\mathfrak c_4}{n^5}
\end{eqnarray*}
\end{small}
\newline
with
\begin{displaymath}
\mathfrak c_4 =
4\mathfrak m_{\mathcal K,\ell}
\mathbb E(\ell(Y_1)^2)
\mathbb E(\exp(\alpha|\ell(Y_1)|)).
\end{displaymath}
So,
\begin{displaymath}
\mathbb E\left(\sup_{K,K'\in\mathcal K_n}
\frac{|U_{K,K',\ell}^{4}(n)|}{n^2}\right)
\leqslant
\frac{\mathfrak c_4}{n^5}.
\end{displaymath}
Therefore,
\begin{displaymath}
\mathbb E\left(\sup_{K,K'\in\mathcal K_n}
\left\{\frac{|U_{K,K',\ell}(n)|}{n^2}
-\frac{\theta}{n}\overline s_{K',\ell}
\right\}\right)
\leqslant
(2^4 + 5.4\mathfrak c_1)
\mathfrak c_2\frac{\log(n)^5}{\theta n} +
2\mathfrak c_3\frac{\log(n)}{n} +
\frac{\mathfrak c_4}{n^5}.
\end{displaymath}
%


%
\subsubsection{Proof of Lemma \ref{bound_trace_term}}
First, the two following results are used several times in the sequel:
\begin{eqnarray}
 \|s_{K,\ell}\|_{2}^{2} & \leqslant &
 \mathbb E(\ell(Y_1)^2)
 \int_{\mathbb R^d}f(x')\int_{\mathbb R^d}
 K(x',x)^2\lambda_d(dx)\lambda_d(dx')
 \nonumber\\
 \label{bound_trace_term_1}
 & \leqslant &
 \mathbb E(\ell(Y_1)^2)\mathfrak m_{\mathcal K,\ell}n
\end{eqnarray}
and
\begin{eqnarray}
 \mathbb E(V_{K,\ell}(n)) & = &
 \mathbb E(\|K(X_1,.)\ell(Y_1) - s_{K,\ell}\|_{2}^{2})
 \nonumber\\
 \label{bound_trace_term_2}
 & = &
 \mathbb E(\|K(X_1,.)\ell(Y_1)\|_{2}^{2}) +
 \|s_{K,\ell}\|_{2}^{2}
 - 2\int_{\mathbb R^d}s_{K,\ell}(x)\mathbb E(K(X_1,x)\ell(Y_1))\lambda_d(dx)
 =\overline s_{K,\ell} -\|s_{K,\ell}\|_{2}^{2}.
\end{eqnarray}
Consider $\mathfrak m(n) := 2\log(n)/\alpha$ and
\begin{displaymath}
v_{K,\ell}(n) := V_{K,\ell}(n) -\mathbb E(V_{K,\ell}(n))
= v_{K,\ell}^{1}(n) + v_{K,\ell}^{2}(n),
\end{displaymath}
where
\begin{displaymath}
v_{K,\ell}^{j}(n) =
\frac{1}{n}\sum_{i = 1}^{n}
(g_{K,\ell}^{j}(n;X_i,Y_i) -\mathbb E(g_{K,\ell}^{j}(n;X_i,Y_i)))
\textrm{ $;$ }
j = 1,2
\end{displaymath}
with, for every $(x',y)\in E$,
\begin{displaymath}
g_{K,\ell}^{1}(n;x',y) :=
\|K(x',.)\ell(y) - s_{K,\ell}\|_{2}^{2}\mathbf 1_{|\ell(y)|\leqslant\mathfrak m(n)}
\end{displaymath}
and
\begin{displaymath}
g_{K,\ell}^{2}(n;x',y) :=
\|K(x',.)\ell(y) - s_{K,\ell}\|_{2}^{2}\mathbf 1_{|\ell(y)| >\mathfrak m(n)}.
\end{displaymath}
On the one hand, by Bernstein's inequality, for any $\lambda > 0$, with probability larger than $1 - 2e^{-\lambda}$,
\begin{displaymath}
|v_{K,\ell}^{1}(n)|
\leqslant
\sqrt{\frac{2\lambda}{n}\mathfrak v_{K,\ell}(n)}
+\frac{\lambda}{n}\mathfrak c_{K,\ell}(n)
\end{displaymath}
where
\begin{displaymath}
\mathfrak c_{K,\ell}(n) =\frac{\|g_{K,\ell}^{1}(n;.)\|_{\infty}}{3}
\textrm{ and }
\mathfrak v_{K,\ell}(n) =
\mathbb E(g_{K,\ell}^{1}(n;X_1,Y_1)^2).
\end{displaymath}
Moreover,
\begin{eqnarray*}
 \mathfrak c_{K,\ell}(n) & = &
 \frac{1}{3}\sup_{(x',y)\in E}
 \|K(x',.)\ell(y) - s_{K,\ell}\|_{2}^{2}
 \mathbf 1_{|\ell(y)|\leqslant\mathfrak m(n)}\\
 & \leqslant &
 \frac{2}{3}\left(
 \mathfrak m(n)^2\sup_{x'\in\mathbb R^d}\|K(x',.)\|_{2}^{2} +\|s_{K,\ell}\|_{2}^{2}
 \right)
 \leqslant
 \frac{2}{3}(\mathfrak m(n)^2 +\mathbb E(\ell(Y_1)^2))\mathfrak m_{\mathcal K,\ell}n
\end{eqnarray*}
by Inequality (\ref{bound_trace_term_1}), and
\begin{eqnarray*}
 \mathfrak v_{K,\ell}(n) & \leqslant &
 \|g_{K,\ell}^{1}(n;.)\|_{\infty}\mathbb E(V_{K,\ell}(n))\\
 & \leqslant &
 2(\mathfrak m(n)^2 +\mathbb E(\ell(Y_1)^2))\mathfrak m_{\mathcal K,\ell}n
 (\overline s_{K,\ell} -\|s_{K,\ell}\|_{2}^{2})
\end{eqnarray*}
by Inequality (\ref{bound_trace_term_1}) and Equality (\ref{bound_trace_term_2}). Then, for any $\theta\in ]0,1[$,
\begin{eqnarray*}
 |v_{K,\ell}^{1}(n)|
 & \leqslant &
 2\sqrt{\lambda
 (\mathfrak m(n)^2 +\mathbb E(\ell(Y_1)^2))\mathfrak m_{\mathcal K,\ell}
 (\overline s_{K,\ell} -\|s_{K,\ell}\|_{2}^{2})}
 +\frac{2\lambda}{3}
 (\mathfrak m(n)^2 +\mathbb E(\ell(Y_1)^2))\mathfrak m_{\mathcal K,\ell}\\
 & \leqslant &
 \theta\overline s_{K,\ell} +
 \frac{5\lambda}{3\theta}
 (1 +\mathbb E(\ell(Y_1)^2))\mathfrak m_{\mathcal K,\ell}\mathfrak m(n)^2
\end{eqnarray*}
with probability larger than $1 - 2e^{-\lambda}$. So, with probability larger than $1 - 2|\mathcal K_n|e^{-\lambda}$,
\begin{displaymath}
S_{\mathcal K,\ell}(n,\theta) :=
\sup_{K\in\mathcal K_n}
\left\{\frac{|v_{K,\ell}^{1}(n)|}{n}
-\frac{\theta}{n}\overline s_{K,\ell}
\right\}
\leqslant
\frac{5\lambda}{3\theta n}
(1 +\mathbb E(\ell(Y_1)^2))\mathfrak m_{\mathcal K,\ell}\mathfrak m(n)^2.
\end{displaymath}
For every $t\in\mathbb R_+$, consider
\begin{displaymath}
\lambda_{\mathcal K,\ell}(n,\theta,t) :=
\frac{t}{\mathfrak m_{\mathcal K,\ell}(n,\theta)}
\textrm{ with }
\mathfrak m_{\mathcal K,\ell}(n,\theta) =
\frac{5}{3\theta n}
(1 +\mathbb E(\ell(Y_1)^2))\mathfrak m_{\mathcal K,\ell}\mathfrak m(n)^2.
\end{displaymath}
Then, for any $T > 0$,
\begin{eqnarray}
 \mathbb E(S_{\mathcal K,\ell}(n,\theta))
 & \leqslant &
 T +\int_{T}^{\infty}\mathbb P(S_{\mathcal K,\ell}(n,\theta)\geqslant
 \lambda_{\mathcal K,\ell}(n,\theta,t)\mathfrak m_{\mathcal K,\ell}(n,\theta))dt
 \nonumber\\
 & \leqslant &
 T + 2|\mathcal K_n|\int_{T}^{\infty}\exp(-\lambda_{\mathcal K,\ell}(n,\theta,t))dt
 \nonumber\\
 & = &
 T + 2|\mathcal K_n|\int_{T}^{\infty}
 \exp\left(-\frac{t}{2\mathfrak m_{\mathcal K,\ell}(n,\theta)}\right)
 \exp\left(-\frac{t}{2\mathfrak m_{\mathcal K,\ell}(n,\theta)}\right)dt
 \nonumber\\
 \label{bound_trace_term_3}
 & \leqslant &
 T + 2\mathfrak c_1|\mathcal K_n|\mathfrak m_{\mathcal K,\ell}(n,\theta)
 \exp\left(-\frac{T}{2\mathfrak m_{\mathcal K,\ell}(n,\theta)}\right)
 \textrm{ with }
 \mathfrak c_1 =
 \int_{0}^{\infty}e^{-r/2}dr = 2.
\end{eqnarray}
Moreover,
\begin{displaymath}
\mathfrak m_{\mathcal K,\ell}(n,\theta)
\leqslant\mathfrak c_2\frac{\log(n)^2}{\theta n}
\textrm{ with }
\mathfrak c_2 =
\frac{20}{3\alpha^2}(1 +\mathbb E(\ell(Y_1)^2))\mathfrak m_{\mathcal K,\ell}.
\end{displaymath}
So, by taking
\begin{displaymath}
T = 2\mathfrak c_2\frac{\log(n)^3}{\theta n},
\end{displaymath}
and since $|\mathcal K_n|\leqslant n$,
\begin{displaymath}
\mathbb E(S_{\mathcal K,\ell}(n,\theta))
\leqslant
2\mathfrak c_2\frac{\log(n)^3}{\theta n} +
4\mathfrak m_{\mathcal K,\ell}(n,\theta)
\frac{|\mathcal K_n|}{n}
\leqslant
6\mathfrak c_2\frac{\log(n)^3}{\theta n}.
\end{displaymath}
On the other hand, by Inequality (\ref{bound_trace_term_1}) and Markov's inequality,
\begin{eqnarray*}
 \mathbb E\left[\sup_{K\in\mathcal K_n}\frac{|v_{K,\ell}^{2}(n)|}{n}\right]
 & \leqslant &
 \frac{2}{n}\mathbb E\left(
 \sup_{K\in\mathcal K_n}\|K(X_1,.)\ell(Y_1) - s_{K,\ell}\|_{2}^{2}\mathbf 1_{|\ell(Y_1)| >\mathfrak m(n)}\right)\\
 & \leqslant &
 \frac{4}{n}\mathbb E\left[\left|
 \ell(Y_1)^2
 \sup_{K\in\mathcal K_n}\|K(X_1,.)\|_{2}^{2}
 +\sup_{K\in\mathcal K_n}\|s_{K,\ell}\|_{2}^{2}\right|^2\right]^{1/2}
 \mathbb P(|\ell(Y_1)| >\mathfrak m(n))^{1/2}
 \leqslant
 \frac{\mathfrak c_3}{n}
\end{eqnarray*}
with
\begin{displaymath}
\mathfrak c_3 = 8\mathfrak m_{\mathcal K,\ell}\mathbb E(\ell(Y_1)^4)^{1/2}\mathbb E(\exp(\alpha|\ell(Y_1)|))^{1/2}.
\end{displaymath}
Therefore,
\begin{displaymath}
\mathbb E\left(\sup_{K\in\mathcal K_n}\left\{
\frac{|v_{K,\ell}(n)|}{n} -\frac{\theta}{n}\overline s_{K,\ell})\right\}\right)
\leqslant 6\mathfrak c_2\frac{\log(n)^3}{\theta n} +\frac{\mathfrak c_3}{n}
\end{displaymath}
and, by Equality (\ref{bound_trace_term_2}), the definition of $v_{K,\ell}(n)$ and Assumption \ref{assumption_K}.(2),
\begin{displaymath}
\mathbb E\left(\sup_{K\in\mathcal K_n}\left\{
\frac{1}{n}|V_{K,\ell}(n) -\overline s_{K,\ell}| -\frac{\theta}{n}\overline s_{K,\ell}\right\}\right)
\leqslant 6\mathfrak c_2\frac{\log(n)^3}{\theta n} +\frac{\mathfrak c_3 +\mathfrak m_{\mathcal K,\ell}}{n}.
\end{displaymath}
%


%
\begin{remark}\label{remark_bound_trace_term}
As mentioned in Remark \ref{remark_risk_bound_main_estimator_2}, replacing the exponential moment condition by the weaker $q$-th moment condition with $q = (12 - 4\varepsilon)/\beta$, $\varepsilon\in ]0,1[$ and $0 <\beta < \varepsilon/2$, allows to get a rate of convergence of order $1/n^{1-\varepsilon}$. Indeed, by Inequality (\ref{bound_trace_term_3}), with $\mathfrak m(n) = n^{\beta}$ and
\begin{displaymath}
T =\frac{2\mathfrak c_1}{\theta n^{1 -\varepsilon}}
\textrm{ with }
\mathfrak c_1 =\frac{5}{3}(1 +\mathbb E(\ell(Y_1)^2))\mathfrak m_{\mathcal K,\ell},
\end{displaymath}
and by letting $\alpha=1+2\beta-\varepsilon$, there exist $n_{\varepsilon,\alpha}\in\mathbb N^*$ and $\mathfrak c_{\varepsilon,\alpha}>0$ not depending on $n$, such that for any $n\geqslant n_{\varepsilon,\alpha}$,
\begin{eqnarray*}
\mathbb E(S_{\mathcal K,\ell}(n,\theta))
&\leqslant &
\frac{2\mathfrak c_1}{\theta n^{1 -\varepsilon}} +4\mathfrak c_1|\mathcal K_n|\frac{n^{2\beta - 1}}{\theta}\exp(-n^{\varepsilon - 2\beta})
\\
& \leqslant &
\frac{2\mathfrak c_1}{\theta n^{1 -\varepsilon}} +4\mathfrak c_1\mathfrak c_{\varepsilon,\alpha}\frac{n^{2\beta}}{\theta n^{\alpha}} =
\frac{2\mathfrak c_1(1+2\mathfrak c_{\varepsilon,\alpha})}{\theta n^{1 -\varepsilon}}.
\end{eqnarray*}
Furthermore, by Markov's inequality,
\begin{displaymath}
\mathbb P(|\ell(Y_1)| > n^{\beta})\leqslant \frac{\mathbb E(|\ell(Y_1)|^{(12 - 4\varepsilon)/\beta})}{n^{12 - 4\varepsilon}}.
\end{displaymath}
So, as previously, there exists a deterministic constant $\mathfrak c_2 > 0$ such that
\begin{displaymath}
\mathbb E\left(\sup_{K,K'\in\mathcal K_n}|W_{K,K',\ell}^{2}(n)|\right)
\leqslant\mathfrak c_2|\mathcal K_n|^2\mathbb P(|\ell(Y_1)| >\mathfrak m(n))^{1/4}
\leqslant\frac{\mathfrak c_3\mathbb E(|\ell(Y_1)|^{(12 - 4\varepsilon)/\beta})^{1/4}}{n^{1-\varepsilon}},
\end{displaymath}
and then
\begin{displaymath}
\mathbb E\left(\sup_{K,K'\in\mathcal K_n}\{
|W_{K,K',\ell}(n)| -\theta\|s_{K',\ell} - s\|_{2}^{2}\}\right)
\leqslant
\frac{\mathfrak c_3}{\theta n^{1-\varepsilon}}
\textrm{ with }
\mathfrak c_3 = 2\mathfrak c_1(1+2\mathfrak c_{\varepsilon,\alpha})+\mathfrak c_2\mathbb E(|\ell(Y_1)|^{(12 - 4\varepsilon)/\beta})^{1/4}.
\end{displaymath}
\end{remark}
%


%
\subsubsection{Proof of Lemma \ref{bound_crossed_term}}
Consider $\mathfrak m(n) = 12\log(n)/\alpha$. For any $K,K'\in\mathcal K_n$,
\begin{displaymath}
W_{K,K',\ell}(n) = W_{K,K',\ell}^{1}(n) + W_{K,K',\ell}^{2}(n)
\end{displaymath}
where
\begin{displaymath}
W_{K,K',\ell}^{j}(n) :=
\frac{1}{n}\sum_{i = 1}^{n}
(g_{K,K',\ell}^{j}(n;X_i,Y_i) -\mathbb E(g_{K,K',\ell}^{j}(n;X_i,Y_i)))
\textrm{ $;$ }
j = 1,2
\end{displaymath}
with, for every $(x',y)\in E$,
\begin{displaymath}
g_{K,K',\ell}^{1}(n;x',y) :=
\langle K(x',.)\ell(y),s_{K',\ell} - s\rangle_2
\mathbf 1_{|\ell(y)|\leqslant\mathfrak m(n)}
\end{displaymath}
and
\begin{displaymath}
g_{K,K',\ell}^{2}(n;x',y) :=
\langle K(x',.)\ell(y),s_{K',\ell} - s\rangle_2
\mathbf 1_{|\ell(y)| >\mathfrak m(n)}.
\end{displaymath}
On the one hand, by Bernstein's inequality, for any $\lambda > 0$, with probability larger than $1 - 2e^{-\lambda}$,
\begin{displaymath}
|W_{K,K',\ell}^{1}(n)|
\leqslant
\sqrt{\frac{2\lambda}{n}\mathfrak v_{K,K',\ell}(n)}
+\frac{\lambda}{n}\mathfrak c_{K,K',\ell}(n)
\end{displaymath}
where
\begin{displaymath}
\mathfrak c_{K,K',\ell}(n) =\frac{\|g_{K,K',\ell}^{1}(n;.)\|_{\infty}}{3}
\textrm{ and }
\mathfrak v_{K,K',\ell}(n) =
\mathbb E(g_{K,K',\ell}^{1}(n;X_1,Y_1)^2).
\end{displaymath}
Moreover,
\begin{eqnarray*}
 \mathfrak c_{K,K',\ell}(n) & = &
 \frac{1}{3}\sup_{(x',y)\in E}
 |\langle K(x',.)\ell(y),s_{K',\ell} - s\rangle_2|
 \mathbf 1_{|\ell(y)|\leqslant\mathfrak m(n)}\\
 & \leqslant &
 \frac{1}{3}\mathfrak m(n)\|s_{K',\ell} - s\|_2
 \sup_{x'\in\mathbb R^d}\|K(x',.)\|_2
 \leqslant
 \frac{1}{3}\mathfrak m_{\mathcal K,\ell}^{1/2}n^{1/2}\mathfrak m(n)\|s_{K',\ell} - s\|_2
\end{eqnarray*}
by Assumption \ref{assumption_K}.(1), and
\begin{displaymath}
\mathfrak v_{K,\ell}(n)\leqslant
\mathbb E(\langle K(X_1,.)\ell(Y_1),s_{K',\ell} - s\rangle_{2}^{2}
\mathbf 1_{|\ell(Y_1)|\leqslant\mathfrak m(n)})\\
\leqslant
\mathfrak m(n)^2
\mathfrak m_{\mathcal K,\ell}\|s_{K',\ell} - s\|_{2}^{2}
\end{displaymath}
by Assumption \ref{assumption_K}.(4). Then, since $\lambda > 0$, for any $\theta\in ]0,1[$,
\begin{eqnarray*}
 |W_{K,K',\ell}^{1}(n)|
 & \leqslant &
 \sqrt{\frac{2\lambda}{n}
 \mathfrak m(n)^2
 \mathfrak m_{\mathcal K,\ell}\|s_{K',\ell} - s\|_{2}^{2}}
 +\frac{\lambda}{3n^{1/2}}\mathfrak m_{\mathcal K,\ell}^{1/2}\mathfrak m(n)\|s_{K',\ell} - s\|_2\\
 & \leqslant &
 \theta\|s_{K',\ell} - s\|_{2}^{2} +
 \frac{\mathfrak m_{\mathcal K,\ell}}{2\theta n}\mathfrak m(n)^2(1 +\lambda)^2
\end{eqnarray*}
with probability larger than $1 - 2e^{-\lambda}$. So, with probability larger than $1 - 2|\mathcal K_n|^2e^{-\lambda}$,
\begin{displaymath}
S_{\mathcal K,\ell}(n,\theta) :=
\sup_{K,K'\in\mathcal K_n}
\{|W_{K,K',\ell}^{1}(n)|
-\theta\|s_{K',\ell} - s\|_{2}^{2}\}
\leqslant
\frac{\mathfrak m_{\mathcal K,\ell}}{2\theta n}\mathfrak m(n)^2(1 +\lambda)^2.
\end{displaymath}
For every $t\in\mathbb R_+$, consider
\begin{displaymath}
\lambda_{\mathcal K,\ell}(n,\theta,t) :=
-1 +
\left(\frac{t}{\mathfrak m_{\mathcal K,\ell}(n,\theta)}\right)^{1/2}
\textrm{ with }
\mathfrak m_{\mathcal K,\ell}(n,\theta) =
\frac{\mathfrak m_{\mathcal K,\ell}}{2\theta n}\mathfrak m(n)^2.
\end{displaymath}
Then, for any $T > 0$,
\begin{eqnarray*}
 \mathbb E(S_{\mathcal K,\ell}(n,\theta))
 & \leqslant &
 T +\int_{T}^{\infty}\mathbb P(S_{\mathcal K,\ell}(n,\theta)\geqslant
 (1 +\lambda_{\mathcal K,\ell}(n,\theta,t))^2\mathfrak m_{\mathcal K,\ell}(n,\theta))dt\\
 & \leqslant &
 T + 2|\mathcal K_n|^2\int_{T}^{\infty}\exp(-\lambda_{\mathcal K,\ell}(n,\theta,t))dt\\
 & = &
 T + 2|\mathcal K_n|^2\int_{T}^{\infty}
 \exp\left(-\frac{t^{1/2}}{2\mathfrak m_{\mathcal K,\ell}(n,\theta)^{1/2}}\right)
 \exp\left(1 -\frac{t^{1/2}}{2\mathfrak m_{\mathcal K,\ell}(n,\theta)^{1/2}}\right)dt\\
 & \leqslant &
 T + 2\mathfrak c_1|\mathcal K_n|^2\mathfrak m_{\mathcal K,\ell}(n,\theta)
 \exp\left(-\frac{T^{1/2}}{2\mathfrak m_{\mathcal K,\ell}(n,\theta)^{1/2}}\right)
 \textrm{ with }
 \mathfrak c_1 =
 \int_{0}^{\infty}e^{1 - r^{1/2}/2}dr.
\end{eqnarray*}
Moreover,
\begin{displaymath}
\mathfrak m_{\mathcal K,\ell}(n,\theta)
\leqslant\mathfrak c_2\frac{\log(n)^2}{\theta n}
\textrm{ with }
\mathfrak c_2 =
\frac{12^2}{2\alpha^2}\mathfrak m_{\mathcal K,\ell}.
\end{displaymath}
So, by taking
\begin{displaymath}
T = 2^3\mathfrak c_2\frac{\log(n)^4}{\theta n},
\end{displaymath}
and since $|\mathcal K_n|\leqslant n$,
\begin{displaymath}
\mathbb E(S_{\mathcal K,\ell}(n,\theta))
\leqslant
2^3\mathfrak c_2\frac{\log(n)^4}{\theta n} +
2\mathfrak c_1\mathfrak m_{\mathcal K,\ell}(n,\theta)
\frac{|\mathcal K_n|^2}{n^2}
\leqslant
(2^3 + 2\mathfrak c_1)\mathfrak c_2\frac{\log(n)^4}{\theta n}.
\end{displaymath}
On the other hand, by Assumption \ref{assumption_K}.(2,4), Cauchy-Schwarz's inequality and Markov's inequality,
\begin{eqnarray*}
 \mathbb E\left(\sup_{K,K'\in\mathcal K_n}|W_{K,K',\ell}^{2}(n)|\right)
 & \leqslant &
 2\mathbb E(\ell(Y_1)^2\mathbf 1_{|\ell(Y_1)| >\mathfrak m(n)})^{1/2}
 \sum_{K,K'\in\mathcal K_n}
 \mathbb E(
 \langle K(X_1,.),s_{K',\ell} - s\rangle_{2}^{2})^{1/2}\\
 & \leqslant &
 2\mathfrak m_{\mathcal K,\ell}^{1/2}\|s_{K',\ell} - s\|_2
 \mathbb E(\ell(Y_1)^4)^{1/4}|\mathcal K_n|^2\mathbb P(|\ell(Y_1)| >\mathfrak m(n))^{1/4}
 \leqslant
 \frac{\mathfrak c_3}{n}
\end{eqnarray*}
with
\begin{displaymath}
\mathfrak c_3 = 2\mathfrak m_{\mathcal K,\ell}^{1/2}
(\mathfrak m_{\mathcal K,\ell}^{1/2} +\|s\|_2)
\mathbb E(\ell(Y_1)^4)^{1/4}
\mathbb E(\exp(\alpha|\ell(Y_1)|))^{1/4}.
\end{displaymath}
Therefore,
\begin{displaymath}
\mathbb E\left(\sup_{K,K'\in\mathcal K_n}\{
|W_{K,K',\ell}(n)| -\theta\|s_{K',\ell} - s\|_{2}^{2}\}\right)
\leqslant (2^3 + 2\mathfrak c_1)\mathfrak c_2\frac{\log(n)^4}{\theta n} +\frac{\mathfrak c_3}{n}
\leqslant
\mathfrak c_4\frac{\log(n)^4}{\theta n}
\end{displaymath}
with $\mathfrak c_4 = (2^3 + 2\mathfrak c_1)\mathfrak c_2 +\mathfrak c_3$.
%


%
\subsection{Proof of Proposition \ref{variance_bound_main_estimator}}
For any $K\in\mathcal K_n$,
\begin{equation}\label{variance_bound_main_estimator_1}
\|\widehat s_{K,\ell}(n;.) - s_{K,\ell}\|_{2}^{2} =
\frac{U_{K,\ell}(n)}{n^2} +\frac{V_{K,\ell}(n)}{n}
\end{equation}
with $U_{K,\ell}(n) = U_{K,K,\ell}(n)$ and $V_{K,\ell}(n) = V_{K,K,\ell}(n)$. Then, by Lemmas \ref{bound_U_statistics} and \ref{bound_trace_term},
\begin{displaymath}
\mathbb E\left(\sup_{K\in\mathcal K_n}\left\{
\left|\|\widehat s_{K,\ell}(n;.) - s_{K,\ell}\|_{2}^{2} -\frac{\overline s_{K,\ell}}{n}\right|
-\frac{\theta}{n}\overline s_{K,\ell}\right\}\right)
\leqslant
\mathfrak c_{\ref{variance_bound_main_estimator}}\frac{\log(n)^5}{\theta n}
\end{displaymath}
with $\mathfrak c_{\ref{variance_bound_main_estimator}} =\mathfrak c_{\ref{bound_U_statistics}} +\mathfrak c_{\ref{bound_trace_term}}$.
%


%
\subsection{Proof of Theorem \ref{risk_bound_main_estimator}}
On the one hand, for every $K\in\mathcal K_n$,
\begin{displaymath}
\|\widehat s_{K,\ell}(n;.) - s\|_{2}^{2} - (1 +\theta)\left(\|s_{K,\ell} - s\|_{2}^{2} +\frac{\overline s_{K,\ell}}{n}\right)
\end{displaymath}
can be written
\begin{displaymath}
\|\widehat s_{K,\ell}(n;.) - s_{K,\ell}\|_{2}^{2} - (1 +\theta)\frac{\overline s_{K,\ell}}{n} +
2W_{K,\ell}(n) -\theta\|s_{K,\ell} - s\|_{2}^{2},
\end{displaymath}
where $W_{K,\ell}(n) := W_{K,K,\ell}(n)$ (see (\ref{W_statistic_definition})). Then, by Proposition \ref{variance_bound_main_estimator} and Lemma \ref{bound_crossed_term},
\begin{displaymath}
\mathbb E\left(\sup_{K\in\mathcal K_n}\left\{
 \|\widehat s_{K,\ell}(n;.) - s\|_{2}^{2} - (1 +\theta)\left(\|s_{K,\ell} - s\|_{2}^{2} +\frac{\overline s_{K,\ell}}{n}\right)\right\}\right)
\leqslant
\mathfrak c_{\ref{risk_bound_main_estimator}}\frac{\log(n)^5}{\theta n}
\end{displaymath}
with $\mathfrak c_{\ref{risk_bound_main_estimator}} =\mathfrak c_{\ref{variance_bound_main_estimator}} +\mathfrak c_{\ref{bound_crossed_term}}$. On the other hand, for any $K\in\mathcal K_n$,
\begin{displaymath}
\|s_{K,\ell} - s\|_{2}^{2} =
\|\widehat s_{K,\ell}(n;.) - s\|_{2}^{2}
-\|\widehat s_{K,\ell}(n;.) - s_{K,\ell}\|_{2}^{2}
- W_{K,\ell}(n).
\end{displaymath}
Then,
\begin{displaymath}
(1 -\theta)\left(\|s_{K,\ell} - s\|_{2}^{2} +\frac{\overline s_{K,\ell}}{n}\right)
-\|\widehat s_{K,\ell}(n;.) - s\|_{2}^{2}
\leqslant
|W_{K,\ell}(n)| -\theta\|s_{K,\ell} - s\|_{2}^{2} +
\Lambda_{K,\ell}(n) -\theta\frac{\overline s_{K,\ell}}{n}
\end{displaymath}
where
\begin{displaymath}
\Lambda_{K,\ell}(n) :=
\left|
\|\widehat s_{K,\ell} - s_{K,\ell}\|_{2}^{2} -\frac{\overline s_{K,\ell}}{n}\right|.
\end{displaymath}
By Equalities (\ref{variance_bound_main_estimator_1}) and (\ref{bound_trace_term_2}),
\begin{displaymath}
\Lambda_{K,\ell}(n) =
\left|\frac{U_{K,\ell}(n)}{n^2} +\frac{v_{K,\ell}(n)}{n} -\frac{\|s_{K,\ell}\|_{2}^{2}}{n}\right|
\end{displaymath}
with $U_{K,\ell}(n) = U_{K,K,\ell}(n)$ (see (\ref{U_statistic_definition})). By Lemmas \ref{bound_trace_term} and \ref{bound_U_statistics}, there exists a deterministic constant $\mathfrak c_1 > 0$, not depending $n$ and $\theta$, such that
\begin{displaymath}
\mathbb E\left(\sup_{K\in\mathcal K_n}\left\{
\Lambda_{K,\ell}(n) -\theta\frac{\overline s_{K,\ell}}{n}\right\}\right)
\leqslant
\mathfrak c_1\frac{\log(n)^5}{\theta n}.
\end{displaymath}
By Lemma \ref{bound_crossed_term},
\begin{displaymath}
\mathbb E\left(\sup_{K\in\mathcal K_n}\{
|W_{K,\ell}(n)| -\theta\|s_{K,\ell} - s\|_{2}^{2}
\}\right)
\leqslant
\mathfrak c_{\ref{bound_crossed_term}}\frac{\log(n)^4}{\theta n}.
\end{displaymath}
Therefore,
\begin{displaymath}
\mathbb E\left(\sup_{K\in\mathcal K_n}\left\{
\|s_{K,\ell} - s\|_{2}^{2} +\frac{\overline s_{K,\ell}}{n} -
\frac{1}{1 -\theta}\|\widehat s_{K,\ell}(n;.) - s\|_{2}^{2}\right\}\right)
\leqslant
\overline{\mathfrak c}_{\ref{risk_bound_main_estimator}}\frac{\log(n)^5}{\theta(1 -\theta) n}
\end{displaymath}
with $\overline{\mathfrak c}_{\ref{risk_bound_main_estimator}} =\mathfrak c_{\ref{bound_crossed_term}} +\mathfrak c_1$.
%


%
\subsection{Proof of Theorem \ref{risk_bound_adaptive_estimator}}
	The proof of Theorem \ref{risk_bound_adaptive_estimator} is dissected in three steps.
	\\
	\\
	\noindent
		\textbf{Step 1.} This first step is devoted to provide a suitable decomposition of 
		\begin{displaymath}
		\|\widehat s_{\widehat K,\ell}(n;\cdot) - s\|_2^2.
		\end{displaymath}
		First,
		\begin{displaymath}
		\|\widehat s_{\widehat K,\ell}(n;\cdot) - s\|_2^2 = \|\widehat s_{\widehat K,\ell}(n;\cdot) - \widehat s_{K_0,\ell}(n;\cdot)\|_2^2 + \|\widehat s_{K_0,\ell}(n;\cdot) - s\|_2^2 - 2\langle \widehat s_{K_0,\ell}(n;\cdot) - \widehat s_{\widehat K,\ell}(n;\cdot) , \widehat s_{K_0,\ell}(n;\cdot) - s\rangle_2 
		\end{displaymath}
		From \eqref{penalty_proposal}, it follows
		that for any $K\in\mathcal K_n$, 
		\begin{eqnarray}\label{oracle_1}
		\nonumber
		\|\widehat s_{\widehat K,\ell}(n;\cdot) - s\|_2^2 
		& \leqslant & \|\widehat s_{K,\ell}(n;\cdot) - s\|_2^2 + \textrm{pen}_{\ell}(K) - \textrm{pen}_{\ell}(\widehat K) + \|\widehat s_{K_0,\ell}(n;\cdot) - s\|_2^2 \\
		\nonumber
		&&- 2 \langle \widehat s_{K,\ell}(n;\cdot) - \widehat s_{\widehat K,\ell}(n\;\cdot) , \widehat s_{K_0,\ell}(n;\cdot) - s \rangle_2\\
		&=& \|\widehat s_{K,\ell}(n;\cdot) - s\|_2^2 +\psi_n(K) -\psi_n(\widehat K)
		\end{eqnarray}
		where 
		\begin{displaymath}
		\psi_n(K) := 2 \langle \widehat s_{K,\ell}(n;\cdot)-s , \widehat s_{K_0,\ell}(n;\cdot) - s \rangle_2- \textrm{pen}_{\ell}(K).\vspace{10pt}\\
		\end{displaymath}
		Let's complete the decomposition of $\|\widehat s_{\widehat K,\ell}(n;\cdot) - s\|_2^2$ by writing
		\begin{displaymath}
		\psi_n(K) = 2(\psi_{1,n}(K) + \psi_{2,n}(K) + \psi_{3,n}(K)),
		\end{displaymath}
		where
		\begin{eqnarray*}
			\psi_{1,n}(K) &:=& \dfrac{U_{K,K_0,\ell}(n)}{n^2},\\
			\psi_{2,n}(K) &:=& -\dfrac{1}{n^2}\left(\displaystyle\sum_{i=1}^n\ell(Y_i)\langle K_0(X_i,.),s_{K,\ell}\rangle_2+ \sum_{i=1}^n\ell(Y_i)
			\langle K(X_i,.), s_{K_0,\ell}\rangle_2 \right)+\dfrac{1}{n}\langle s_{K_0,\ell},s_{K,\ell} \rangle_2
			\textrm{ and}\\
			\psi_{3,n}(K) &:=& W_{K,K_0,\ell}(n)+ W_{K_0,K,\ell}(n) + \langle  s_{K,\ell}- s, s_{K_0,\ell}- s \rangle_2.\\
		\end{eqnarray*}	
		\textbf{Step 2.} In this step, we give controls of the quantities
		\begin{displaymath}
		\mathbb E(\psi_{i,n}(K)) \; \text{and} \; \mathbb E(\psi_{i,n}(\widehat K)) \; ; \; i=1,2,3. 
		\end{displaymath}
		\begin{itemize}
			\item %
			By Lemma \ref{bound_U_statistics}, for any $\theta\in ]0,1[$,
			\begin{equation*}
			\mathbb E(|\psi_{1,n}(K)|)
			\leqslant \frac{\theta}{n}\overline s_{K,\ell} +
			\mathfrak c_{\ref{bound_U_statistics}}\frac{\log(n)^5}{\theta n}
			\end{equation*}
			and
			\begin{equation*}
			\mathbb E(|\psi_{1,n}(\widehat K)|)
			\leqslant \frac{\theta}{n}\mathbb E(\overline s_{\widehat K,\ell}) +
			\mathfrak c_{\ref{bound_U_statistics}}\frac{\log(n)^5}{\theta n}. 
			\end{equation*}
			\item On the one hand, for any $K,K'\in\mathcal K_n$, consider 
			\begin{displaymath}
			\Psi_{2,n}(K,K') := \dfrac{1}{n}\displaystyle\sum_{i=1}^n\ell(Y_i)\langle K(X_i,.),s_{K',\ell}\rangle_2.
			\end{displaymath}
			Then, by Assumption \ref{additional_assumption_K},
			\begin{eqnarray*}
				\mathbb E\left(\sup_{K,K'\in\mathcal K_n}|\Psi_{2,n}(K,K')|\right) 
				&\leqslant&
				\mathbb E(\ell(Y_1)^2)^{1/2}
				\mathbb E\left(\sup_{K,K'\in\mathcal K_n}
				\langle K(X_1,.),s_{K',\ell}\rangle_{2}^{2}\right)^{1/2}\\
				&\leqslant&
				\overline{\mathfrak m}_{\mathcal K,\ell}^{1/2}\mathbb E(\ell(Y_1)^2)^{1/2}.
			\end{eqnarray*}
			On the other hand, by Assumption \ref{assumption_K}.(2),
			\begin{displaymath}
				|\langle s_{K,\ell},s_{K_0,\ell} \rangle_2| 
				\leqslant  \mathfrak{m}_{\mathcal K,\ell}.
			\end{displaymath}	
			Then, there exists a deterministic constant $\mathfrak c_1 > 0$, not depending on $n$ and $K$, such that
			\begin{equation*}
			\mathbb E(|\psi_{2,n}(K)|)\leqslant \frac{\mathfrak c_1}{n}
			\textrm{ and }
			\mathbb E(|\psi_{2,n}(\widehat K)|)\leqslant \frac{\mathfrak c_1}{n}. 
			\end{equation*}
			\item By Lemma \ref{bound_crossed_term},
			\begin{eqnarray*}
				\mathbb E(|\psi_{3,n}(K)|)
				&\leqslant & \dfrac{\theta}{4}(\|s_{K,\ell} - s\|_{2}^{2}+\|s_{K_0,\ell} - s\|_{2}^{2})+8\mathfrak c_{\ref{bound_crossed_term}}\frac{\log(n)^4}{\theta n}\\
				&&+\left(\dfrac{\theta}{2}\right)^{1/2}\|s_{K,\ell} - s\|_{2}\times  \left( \dfrac{2}{\theta}\right)^{1/2} \|s_{K_0,\ell} - s\|_{2}\\
				&\leqslant & \dfrac{\theta}{2}\|s_{K,\ell} - s\|_{2}^{2} + \left(\dfrac{\theta}{4}+\dfrac{1}{\theta}\right)\|s_{K_0,\ell} - s\|_{2}^{2}+8\mathfrak c_{\ref{bound_crossed_term}}\frac{\log(n)^4}{\theta n}
			\end{eqnarray*}	
			and
			\begin{equation*}
			\mathbb E(|\psi_{3,n}(\widehat K)|)\leqslant \frac{\theta}{2}\mathbb E(\|s_{\widehat K,\ell} - s\|_{2}^{2}) + \left(\dfrac{\theta}{4}+\dfrac{1}{\theta}\right)\|s_{K_0,\ell} - s\|_{2}^{2}+8\mathfrak c_{\ref{bound_crossed_term}}\frac{\log(n)^4}{\theta n}.\vspace{3pt}\\
			\end{equation*}
		\end{itemize}
		\noindent
		\textbf{Step 3.}
		\noindent
		By the previous step, there exists a deterministic constant $\mathfrak c_2 > 0$, not depending on $n$, $\theta$, $K$ and $K_0$, such that
		\begin{displaymath}
			\mathbb E(|\psi_n(K)|)
			\leqslant \theta \left(\|s_{K,\ell} - s\|_{2}^{2} +\dfrac{\overline s_{K,\ell}}{n}\right) +\left(\dfrac{\theta}{2}+\dfrac{2}{\theta}\right)\|s_{K_0,\ell} - s\|_{2}^{2} +\mathfrak c_2\dfrac{\log(n)^5}{\theta n}
		\end{displaymath} 
		and
		\begin{displaymath}
		\mathbb E(|\psi_n(\widehat K)|)
		\leqslant \theta\mathbb E\left(\|s_{\widehat K,\ell} - s\|_{2}^{2} +\dfrac{\overline s_{\widehat K,\ell}}{n}\right) +\left(\dfrac{\theta}{2}+\dfrac{2}{\theta}\right)\|s_{K_0,\ell} - s\|_{2}^{2} +\mathfrak c_2\dfrac{\log(n)^5}{\theta n}.
		\end{displaymath}  
		Then, by Theorem \ref{risk_bound_main_estimator},
		\begin{eqnarray*}
			\mathbb E(|\psi_n(K)|)
			&\leqslant& \dfrac{\theta}{1-\theta}
			\mathbb E(\|\widehat s_{K,\ell}(n;.) - s\|_{2}^{2}) +\left(\dfrac{\theta}{2}+\dfrac{2}{\theta}\right)\|s_{K_0,\ell} - s\|_{2}^{2}+\left(\dfrac{\mathfrak c_2}{\theta} +\dfrac{\mathfrak c_{\ref{risk_bound_main_estimator}}}{1 -\theta}\right)\dfrac{\log(n)^5}{n}
		\end{eqnarray*}	
		and
		\begin{eqnarray*}
			\mathbb E(|\psi_n(\widehat K)|)
			&\leqslant& \dfrac{\theta}{1-\theta}
			\mathbb E(\|\widehat s_{\widehat K,\ell}(n;.) - s\|_{2}^{2}) +\left(\dfrac{\theta}{2}+\dfrac{2}{\theta}\right)\|s_{K_0,\ell} - s\|_{2}^{2}+\left(\dfrac{\mathfrak c_2}{\theta} +\dfrac{\mathfrak c_{\ref{risk_bound_main_estimator}}}{1 -\theta}\right)\dfrac{\log(n)^5}{n}.
		\end{eqnarray*}	
		By decomposition \eqref{oracle_1}, there exist two deterministic constants $\mathfrak c_3,\mathfrak c_4 >0$, not depending on $n$, $\theta$, $K$ and $K_0$, such that 
		\begin{eqnarray*}
			\mathbb E(\|\widehat s_{\widehat K,\ell}(n;\cdot)-s\|_2^2) 
			&\leqslant & \mathbb E(\|\widehat s_{K,\ell}(n;\cdot)-s\|_2^2) +\mathbb E(| \psi_n(K)|) + \mathbb E(|\psi_n(\widehat K)|)\\
			&\leqslant & \left(1 + \dfrac{\theta}{1-\theta}\right)	\mathbb E(\|\widehat s_{K,\ell}(n;\cdot)-s\|_2^2) +  \dfrac{\theta}{1-\theta}\mathbb E(\|\widehat s_{\widehat K,\ell}(n;.) - s\|_{2}^{2})\\
			&&+\dfrac{\mathfrak c_3}{\theta}\|s_{K_0,\ell} - s\|_{2}^{2}+\dfrac{\mathfrak c_4}{\theta(1-\theta)}\cdot\dfrac{\log(n)^5}{n}.
		\end{eqnarray*}	
		This concludes the proof.
\\
\\
\textbf{Acknowledgments.} This project has received funding from the European Union's Horizon 2020 research and innovation programme under grant agreement N$^{\circ}$811017. The authors want also to thank Fabienne Comte for her careful reading and advices.
%


%

%
\end{document}